\documentclass[opre,copyedit]{informs1}

\usepackage{natbib}
 \bibpunct[, ]{(}{)}{,}{a}{}{,}%
 \def\bibfont{\small}%
 \def\bibsep{\smallskipamount}%
 \def\bibhang{24pt}%
\usepackage{hyperref}

\usepackage{appendix}
\usepackage{times,amsmath,amssymb,psfrag,subfigure}
\usepackage{graphicx,epsfig} 
\usepackage{graphics,color} 

\numberwithin{equation}{section}

\newtheorem{lemma}{Lemma}[section]

\newtheorem{assumption}{Assumption}[section]

\def\v0{\vec{0}}
\def\<{\left<}
\def\>{\right>}
\def\({\left(}
\def\){\right)}

\def\bB{\mathbf{B}}

\def\bD{\mathbf{\Delta}}

\def\bI{\mathbf{I}}

\def\cC{\mathcal{C}}

\def\cK{\mathcal{K}}

\def\cP{\mathcal{P}}
\def\cQ{\mathcal{Q}}
\def\cR{\mathcal{P}}
\def\cS{\mathcal{S}}

\def\cV{\mathcal{V}}

\def\reals{\mathbb{R}}

\def\ints{\mathbb{Z}}

\newtheorem{Theorem}{Theorem}[section]

\newcommand{\eq}{\begin{equation}}
\newcommand{\en}{\end{equation}}
\newcommand{\eqm}{\begin{multline}}
\newcommand{\enm}{\end{multline}}
\newtheorem{example}{Example}

\newcommand{\bTheorem}{ \begin{Theorem}  }
\newcommand{\eTheorem}{ \end{Theorem}    }
\newtheorem{Proposition}{Proposition}[section]
\newcommand{\bProposition}{ \begin{Proposition}  }
\newcommand{\eProposition}{ \end{Proposition}    }
\newtheorem{Lemma}{Lemma}[section]
\newcommand{\bLemma}{ \begin{Lemma}  }
\newcommand{\eLemma}{ \end{Lemma}    }
\newtheorem{Corollary}{Corollary}[section]
\newcommand{\bCorollary}{ \begin{Corollary}  }
\newcommand{\eCorollary}{ \end{Corollary}    }
\newtheorem{Remark}{Remark}[section]
\newcommand{\bRemark}{ \begin{Remark} \rm }
\newcommand{\eRemark}{ \end{Remark}    }
\newtheorem{Definition}{Definition}[section]
\newcommand{\bDefinition}{ \begin{Definition} \rm }
\newcommand{\eDefinition}{ \end{Definition}    }
\newcommand{\bProof}{ \noindent {\sc Proof:} }
\newcommand{\eProof}{\hfill $\blacksquare$ }
\newcommand{\bProofof}[1]{ {\sc Proof of #1:} \hskip 5pt}

\newcommand{\go}{{\rightarrow}}

\begin{document}
\ARTICLEAUTHORS{%
\AUTHOR{Carri W. Chan}
\AFF{Division of Decision, Risk and Operations, Columbia Business School \EMAIL{cwchan@columbia.edu}}
\AUTHOR{Mor Armony}
\AFF{Stern School of Business, New York University \EMAIL{marmony@stern.nyu.edu}}
\AUTHOR{Nicholas Bambos}
\AFF{Departments of Electrical Engineering and Management Science \& Engineering, Stanford University \EMAIL{bambos@stanford.edu}}
} 


\TITLE{Fairness in overloaded parallel queues}
\ABSTRACT{
Maximizing throughput for heterogeneous parallel server queues has received quite a bit of attention from the research community and the stability region for such systems is well understood.  However, many real-world systems have periods where they are temporarily overloaded.  Under such scenarios, the unstable queues often starve limited resources.  This work examines what happens during periods of temporary overload.  Specifically, we look at how to fairly distribute stress.  We explore the dynamics of the queue workloads under the MaxWeight scheduling policy during long periods of stress and discuss how to tune this policy in order to achieve a target fairness ratio across these workloads.
}

\maketitle
\section{Introduction}
Queueing systems are generally designed with sufficient allocation of resources  to ensure that the system is stable. Subsequently, a large body of research in queueing has focused
on characterizing the stability region and on evaluating and
optimizing  performance when the queues are operated within this
region. However, in reality, it is inevitable that a system will
repeatedly enter long periods of overload, either due to an
unpredictable increase in load or because some service resources may
become unavailable/breakdown. In such periods, which may last for a fairly long duration of time, the system may have
more load than it can handle.
Hence, it is natural to want to understand what happens during these
periods of {\em temporary instability}.

Some examples of real systems which experience periods of
instability can be found in internet communication networks,
hospitals, and call centers. Consider systems where traffic is
bursty.  For instance, in communications networks, with the release
of a new online game, the system may be over-stressed as many people
attempt to simultaneously enter the virtual world.  At a hospital, there may be a surge in demand for the limited hospital resources during a catastrophic event or viral outbreak.  In these cases, the system load may increase
to potentially unstable levels. A separate scenario of temporary overload may occur when resources are lost. In call centers, a power outage may
effectively remove a pool of servers.  Consequently, while the system load may
remain constant, the available service resources have been
reduced, stressing the system.

In long periods of overload, it is inevitable that the
system backlog will grow significantly, so the question of
performance optimization becomes moot. However, it is desirable to
share the stressed resources in a {\em fair} manner among the
different customers. In this paper, we study fair resource
allocation in multiclass parallel server queues during long periods
of overload or temporary instability.

How does one define fairness in an overloaded system? Consider a
system which is overloaded for a very long time window. During this
time window, the system is essentially unstable, so that the total
backlog grows without bound. However, it is plausible that some
customer classes will enjoy a greater allocation of the stressed
resources, enough such that the system will appear stable to these
classes. This is unfair. Our {\em novel} notion of fairness requires
that the backlog of the different customer classes will grow large
according to a predetermined proportion, thereby proportionately distributing stress across customer classes. This is the same as
requiring that the {\em backlog vector} will grow along a
prespecified {\em direction}.

Working with this notion of fairness, our ultimate goal is to
determine how to operate the overloaded system such that the total
backlog is minimized subject to a fairness
constraint. Our main result is that the MaxWeight scheduling
algorithm, with carefully selected weight parameters, will achieve this goal asymptotically, as the duration of the overloaded period grows
large. The MaxWeight policy selects, at each time point, a service
configuration which maximizes the aggregate service rate weighted by
the queue backlogs and a controllable weight vector
\citep{tassiulas_tac92,Mandelbaum_Stolyar}. MaxWeight has been shown
to be throughput maximizing (see for example
\cite{tassiulas_tac95,tassiulas_sm00,mckeown_tc99,armony_questa03}).
That is, as long as the system is stabilizable, it will be stable
under MaxWeight.

Our analysis first shows that, under MaxWeight with arbitrary fixed
weights, the backlog in an overloaded system will grow to infinity
according to a well defined direction. We then proceed to show that this direction can be turned to any feasible fairness proportion by
carefully selecting the weight parameters of MaxWeight. Finally, we
establish that, in the limit, MaxWeight will minimize the total
backlog among all algorithms that obtain this fairness direction.

The approach we take in our analysis is a trace-based approach
\citep{loynes, armony_questa03, ross_ton09} which analyzes individual
traces of customers' arrival times and workload requests. To show
that the backlog grows along a specific direction, we use {\em
direct} geometric arguments that examine the dynamics of an actual
trajectory of the backlogs and how it evolves over time. In
particular, our statistical assumptions are extremely mild, and only
require that the workload entering the system has a well defined
long time average.

The rest of the paper is organized as follows: We conclude this
introduction by surveying some closely related literature. In
Section \ref{sec:model} we formally introduce our queueing model and the notion of fairness which we study.  In Section \ref{sec:asym_dyn} we show the existence of a unique limit for the queue backlogs when the system is unstable.  In Section \ref{sec:control} we discuss how to control the backlog to achieve our desired fairness criterion.  In Section \ref{sec:numerical} we show some numerical results to demonstrate the performance of our proposed algorithm in practice.  Finally, we conclude in Section \ref{sec:conc}.

\noindent \textbf{Related Literature}

\cite{Stolyar_04} analyzes the same type of parallel
server queueing systems, under the assumption that the system is in
heavy-traffic. That is,  the system is stable but operates close
to the boundary of the stability region. For such systems, the
author shows that under MaxWeight and a complete resource pooling
condition, the backlog vector will experience a state-space
collapse, which implies that it will asymptotically follow a well
defined direction. The author also shows that MaxWeight minimizes
the workload, where the latter is defined as the total backlog
weighted by the above direction. These results are similar to our
convergence and backlog minimization results, but since the system
considered in \cite{Stolyar_04} is not overloaded, much stronger
assumptions are required there--such as the Markov property and the
complete resource pooling condition.
\cite{Mandelbaum_Stolyar} generalizes these results for convex
holding cost functions and a corresponding generalized $c\mu$ rule.

More recently, Shah and Wischik in \citep{Shah_Wischik} study the
asymptotic behavior of {\em fluid models} in overload under
MaxWeight and other policies. Analogously to our result, they show
that as time grows large, the fluid vector will approach a well
defined direction. While the network structure and the set of
policies that \cite{Shah_Wischik} studies are more general than our
setting, that paper focuses on {\em characterizing} the asymptotic
direction of the backlog of the overloaded system.  Our paper, on the other hand, focuses on {\em controlling} this direction to satisfy a {\em
fairness} constraint and on minimizing the backlog subject to this constraint. Methodologically the two papers are also
different; \cite{Shah_Wischik} uses Lyapunov functions to prove
convergence of the continuous {\em fluid} trajectory, while our analysis
directly examines the dynamics of the {\em queueing trajectory}.

Also in the networked setting, \cite{georgiadis_transit_06} consider overloaded sensor networks where transmissions are scheduled in a distributed manner.  Specifically, \cite{georgiadis_transit_06} show that a policy analogous to the Adaptive Back Pressure policy of \cite{tassiulas_tac95} (reference [20] in Tassiulas' paper) obtains, in the limit, a backlog vector which is the so-called "most balanced" among all feasible limits, in overload. Similar to \cite{Shah_Wischik}, the authors consider a fluid model rather than the direct {\em queueing trajectory} dynamics considered in this work.

Our concept of fairness in this paper is also somewhat unusual. Two commonly
used fairness criteria are Max-Min and Proportional Fairness.   A
system is considered Max-Min fair if the utility of the user with
the minimum utility is maximized.  On the other hand, a system is
considered to be proportionally fair if the amount of service
resources a user is allocated corresponds to the proportion of
anticipated resource consumption required by the user. These two
notions of fairness can be tied together under the description of
$\alpha$-fairness with different values of $\alpha$ \citep{mo_ton00}.
There has been a substantial amount of work in the development of
scheduling policies which ensure some fairness criteria (see
\cite{mo_ton00,kelly_aap04,eryilmaz_infocom05,neely_infocom05,eryilmaz_jsac06,neely_jsac06,bonald_questa06,massoulie_aap07}
and related works). A common thread in those papers is an assumption of stability, whereas we consider fairness during periods of temporary instability. We define fairness by a set of ratios which specifies the desired proportion of the aggregate backlog each queue contributes.  In Section \ref{ssec:otherfairness} we elaborate more
on the connection between the traditional notions of fairness and
the one used here.

\cite{perry_ms_09,perry_or_11} consider how to deal with unexpected overload.  They consider overload in the setting of two initially separated service systems which share resources when one becomes overloaded.  Using a fluid model, they consider how to share resources across the two systems in order to maintain a constant ratio between the two queues.  Our setting is quite different as there is a  single system with multiple parallel queues.  However, we share a similar control goal of maintaining a constant ratio between queues, of which we may have more than two.

\section{The Queueing Model}\label{sec:model}
We consider a queueing system with $Q$ queues and $N$ service vectors.  New jobs arrive to the system and are queued up to be served.  The system administrator dynamically selects which service configuration to implement for service of available jobs in the queue.  We consider how to make this selection in a `fair' manner, where we will be more precise as to what we mean by fair in the coming discussion.

More formally, we consider a queueing system with $Q$ parallel queues, indexed by $q \in \cQ = \{1,2,\dots,Q\}$.  Time is discrete and indexed by $t\in \ints^+_0$.  Let $A_q(t)$ be the number of jobs arriving to queue $q$ in time slot $t$.  We assume that $0 \leq A_q(t) \leq \bar{A}_q$ for some arbitrarily large $\bar{A}_q > 0$.  For each $q \in \cQ$, we assume that
\eq
\lim_{t\rightarrow \infty} \frac{\sum_{s=0}^{t-1} A_q(s)}{t} = \rho_q \in (0,\infty)
\en
is well-defined, positive, and finite.  The arrival vector in time slot $t$ is thus
\eq
A(t) = (A_1(t),A_2(t),\dots,A_q(t),\dots,A_Q(t))
\en
with corresponding long-term traffic load vector:
\eq
\rho = (\rho_1,\rho_2,\dots,\rho_q,\dots,\rho_Q)
\en
Arriving jobs are buffered in their respective queues and are served in a first-come-first-served (FCFS) manner within each queue.

In each time slot, a service vector $S\in \cS = \{S_1,S_2,\dots,S_n,\dots,S_N\}$ is selected to be used.  Each service vector is a $Q$-dimensional vector:
\eq
S= (S_{1},S_{2},\dots,S_{q},\dots,S_{Q})
\en
where $S_{q} \geq 0$ is the number of jobs removed from queue $q$ in a single time slot when service configuration $S$ is used.

The workload in queue $q$ at time $t$ is denoted by $X_q(t)$ with corresponding workload vector:
\eq
X(t) = (X_1(t),X_2(t),\dots,X_q(t),\dots,X_Q(t))
\en
This corresponds to the number of jobs in each queue in time-slot $t$.

Given workload vector $X(t)$ and service configuration $S(t)$, the number of jobs that are served and depart from queue $q$ is:
\eq
D_q(t) = \min \{ S_q(t),X_q(t)\}
\en
where the minimum accounts for the fact that jobs can only be serviced if they are already waiting in the queue.  Hence, if $D_q(t) = X_q(t) < S_q(t)$, there is some idle service provided by service vector $S(t)$ due to the lack of available jobs to be processed.  The workload vector evolves as:
\eq
X(t+1) = X(t) + A(t) - D(t)
\en
Assuming the queue begins empty at time $t = 0$, then the workload vector $X(t)$ is given by:
\eq
X(t) = \sum_{s=0}^{t-1} A(s) -  \sum_{s=0}^{t-1} D(s)
\en

\subsubsection*{Applications}
We note that many applications of interest can be modeled in this way.  For example, in communication networks, the backlogs correspond to packets waiting to be transmitted on various wired or wireless links while the service vectors correspond to various packet switch configurations.  In call centers, the backlogs correspond to the number of customers of various classes waiting to be served while the service vectors correspond to specific allocations of staff with differing skills to each customer class.

\subsection{(In)Stability Region}\label{ssec:stabilityregion}
Loosely speaking, we consider a system to be rate stable if the average job departure rate is equal to the average arrival rate. We define the stability region of a system, $\cP$, such that if $\rho \in \cP$, then there exists some policy which guarantees for each queue, $q$, that:
\eq
\lim_{t \rightarrow \infty} \frac{X_q(t)}{t} = 0
\en
The stability region can be characterized as:
\begin{eqnarray}
\cP &=&\left\{ \rho \in \reals^Q_+ : \rho \leq \sum_{i=1}^N \alpha_i S_i\textrm{ for some }\alpha_i \geq 0, S_i\in \cS\textrm{ such that }\sum_{i=1}^N \alpha_i  =1 \right\} \nonumber \\
&=&\left\{\rho\in\reals_+^Q: \<\rho,\Delta v\> \leq \max_{S\in\cS}\<S,\Delta v\> \mbox{ for every } v\in\reals^Q_+ \right\}
\end{eqnarray}
as defined in \cite{ross_ton09}.  Hence, any system load which is dominated by a convex combination of service vectors is stablizable.  Furthermore, stability is ensured by using this convex combination.

When the system load is outside of the stability region, then the system is unstable.  Hence,  if $\rho \not \in \cP$, then there does not exist any policy which guarantees that:
\eq
\lim_{t \rightarrow \infty} \frac{X_q(t)}{t} = 0
\en
and with probability $1$,
\eq
\lim_{t \rightarrow \infty} \frac{X_q(t)}{t} > 0
\en
for some $q$ (see \cite{armony_questa03}). We refer to this as the \emph{instability} region.  Consider the following simple example of the stability and instability regions.
\begin{example}
Consider an example with two queues ($Q = 2$) and three ($N = 3$) service vectors:
\eq
S_1 = \left(
  \begin{array}{c}
    2 \\
    1 \\
  \end{array}
\right),
S_2 = \left(
  \begin{array}{c}
    1 \\
    1.5 \\
  \end{array}
\right),
S_3 = \left(
  \begin{array}{c}
    1 \\
    1 \\
  \end{array}
\right)
\en
The stability region, $\cP$, is given by the convex hull of $S_1,S_2,S_3$ (and their projections onto the axes).  Note that $S_3$ is a {\em non-essential} service vector, since it is dominated by a convex combination of $S_1$ and $S_2$ and its removal would not change $\cP$.  Any load vector $\rho$ within the shaded region in Figure \ref{fig:twoqueue_example} is stabilizable; outside the region, is not.  Hence, $\hat{\rho}$ is stabilizable while $\rho'$ is not.
\begin{figure}[h]
\psfrag{rh}{$\hat{\rho}$}
\psfrag{rd}{${\rho}'$}
\psfrag{S1}{$S_1$}
\psfrag{S2}{$S_2$}
\psfrag{S3}{$S_3$}
\begin{center}
\includegraphics[scale=.5]{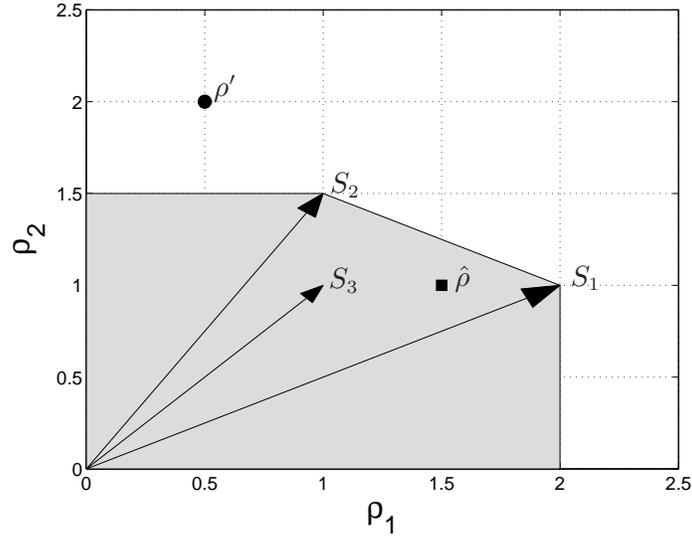}
\caption{An example of the stability region for a 2-queue queueing system.} \label{fig:twoqueue_example}
\end{center}
\end{figure}
\end{example}

The focus of this work is the understand and control the dynamics of systems which are temporarily overloaded over a long time horizon. To do this, we approximate the dynamics of a {\em temporarily} unstable system with a system which is {\em permanently} unstable ($\rho \not\in \cP$).  When the period of overload  is long enough, the dynamics of the two systems will be similar.  Hence, our analysis focuses on the dynamics of unstable systems.

\subsection{Fairness}\label{ssec:fairness}
The goal of this work is to determine a service discipline which allows us to fairly serve queues when the system is temporarily overloaded ($\rho \not \in \cP$).  There are many different definitions of fairness (see for instance \cite{mazumdar_tcomm91,kelly_jors98,mo_ton00}).  In this paper, we focus on a notion of fairness where the workloads grow according to a fixed proportion.  More formally, we assume we are given a set of ratios $\theta_q\geq0, \sum_q \theta_q = 1$.  These ratios specify the proportion of the aggregate workload which each queue contributes.  Whenever $\rho \not \in \cP$,  the goal is to control the workload such that for large  $t$:
\eq
\frac{X_q(t)}{\sum_k X_k(t)} \approx \theta_q, \forall q \in \cQ
\en
 More precisely, we want to control $X(t)$ such that:
\eq
\lim_{t\rightarrow\infty}\frac{X(t)}{t} = \eta
\en
where $\eta$ satisfies our fairness criterion:
\eq
\frac{\eta_i}{\sum_j \eta_j} = \theta_i, \forall i
\en
$\theta$ defines a \emph{direction} which we want the backlogs to grow; $\eta$ includes the scale factor to specify how quickly they grow.

There may be many policies which achieve this fairness criterion and our goal is to find the policy which achieves it with the \textit{minimal workload}.  Let $X(t)$ correspond to the workload at time $t$ given service vector $S(t)$ is used in time slot $t$.  Our goal is to find $S(t)$ such that for large  $t$:
\eq \label{eq:faircrit}
\lim_{t\go\infty}\frac{X_q(t)}{\sum_k X_k(t)} = \theta_q, \forall q \in \cQ
\en
Furthermore, if there exists $X'(t)$ which uses $S'(t)$ in time slot $t$ and $\lim_{t\go\infty}\frac{X'_q(t)}{\sum_k X'_k(t)} = \theta_q, \forall q \in \cQ$, then
\eq \label{eq:min}
\lim_{t\rightarrow \infty} \frac{X(t)}{t} \leq \lim_{t\rightarrow \infty} \frac{X'(t)}{t}
\en

We will go into more detail about this fairness criterion in Section \ref{sec:control}.  Further, we will show that this definition can be generalized to other fairness definitions.  Our proposal is to use MaxWeight scheduling policies to achieve our fairness criterion \cite{armony_questa03,ross_ton09}.

\subsection{The MaxWeight Scheduling Policy}\label{ssec:MWM}
The stability region guarantees the existence of a stabilizing policy.  We now briefly review a family of stabilizing policies, which we refer to as MaxWeight Scheduling.  This is the special case of Projective Cone Scheduling (PCS) when projection matrices are diagonal. This family of policies is characterized by service cones wherein the set of workloads $X(t)$ under which service vector $S_i$ is used forms a cone in the workload-space.  We focus on this family of policies because they work `well', i.e. are stabilizing when possible, they are simple to implement, and as we will see later, they are optimal in the sense that they obtain the minimum backlog \eqref{eq:min} subject to the fairness criterion \eqref{eq:faircrit}.

Let $\bB $ be a symmetric, positive definite matrix with non-positive off-diagonal elements.  Let $X = X(t)$ be the workload at time $t$.  Then, the PCS algorithm selects service vector $S^*$:
\begin{eqnarray}
S^* \in S^{\max}(X) = \arg\max_{S\in \cS} \< S, \bB  X \> = \arg\max_{S \in \cS} X^T\bB  S
\end{eqnarray}
It is shown in \cite{ross_ton09} that the PCS policy is rate stable for any $\rho \in \cP$.  This policy is also desirable from an implementation standpoint as it only requires current information (workload state) and does not require knowledge of the system load $\rho$.

We define the service cone $C_i$ as the set of workload vectors, $X$, such that service vector $S_i$ is used under the  algorithm defined by positive definite matrix $\bB $; i.e.:
\eq
C_i = C_i(\bB ) = \{X|S_i \in \arg\max_{S\in \cS} \< S, \bB  X \>\}
\en
For the rest of the discussion, we will suppress the dependence of the service cones on the matrix $\bB $.

Throughout this work, we will focus  on the case where the  matrix $\bB$   is a diagonal matrix, $\bD$ (and positive-definite, hence, all its diagonal elements are positive).  This specifies the set of algorithms  to MaxWeight scheduling algorithms, which have been studied extensively in the past.
\begin{assumption}\label{as:Bmatrix}
$\bB = \bD $ is a diagonal, positive definite matrix, i.e. $\bD _{qq} > 0$ and $\bD _{qq'} = 0$, $\forall q \not = q'$
\end{assumption}

While we are ultimately interested in the behavior of this queueing system when it is temporarily overloaded, our focus on MaxWeight scheduling policies ensures that the system is stabilized when possible.  We will examine the asymptotic behavior of this class of policies and its implications in term of our fairness criterion.  We will see that this family of algorithms allows us to satisfy certain fairness criterion by manipulating the $\bD $ matrix. Furthermore, we will see in Section \ref{ssec:MWMfair} that the MaxWeight scheduling policies achieves the lowest workload which satisfies our fairness criterion.   In order to see this, we must first understand how the $\bD $ matrix affects the asymptotic behavior of our queueing system.

\section{Asymptotic Dynamics in Overload}\label{sec:asym_dyn}
Before we consider how to control the workload vector, $X(t)$, we first begin by building an understanding of its asymptotic dynamics given diagonal MaxWeight  matrix $\bD$.   We already know that if $\rho \in \cP$, then the MaxWeight scheduling policy stabilizes the system and the workload is always finite.  On the other hand, many real systems enter periods of instability where the system load is temporarily outside of the stability region. From \cite{armony_questa03}, when $\rho \not \in \cP$ the workload explodes: $X(t) \rightarrow \infty$.  However, we are interested in finding out how exactly this happens.  Is there a finite limit for $\lim_{t \rightarrow \infty} \frac{X(t)}{t}$?  If so, what is it and what does it depend on?  This section is devoted to answering these questions.

Suppose that the system operates in overload, in the sense that $\rho\notin\cR$, where
\eq
\cR=\left\{\rho\in\reals_+^Q: \<\rho,\Delta v\> \leq \max_{S\in\cS}\<S,\Delta v\> \mbox{ for every } v\in\reals^Q \right\},
\en
as defined in \cite{ross_ton09}. 
We consider the limit defined below:
\eq \label{eq:def_eta}
H = \limsup_{t\go\infty} \<\frac{X(t)}{t},\bD \frac{X(t)}{t}\>
\en
and select a convergent increasing unbounded subsequence $\{t_c\}$ on which the `$\limsup$' is attained\footnote{See footnote 12 of \cite{ross_ton09} concerning why such a subsequence exists} -- hence,
\eq
\lim_{c\go\infty}\frac{X(t_c)}{t_c} = \eta
\en
and
\eq
\lim_{c\go\infty}\<\frac{X(t_c)}{t_c},\bD \frac{X(t_c)}{t_c}\> = \<\eta,\bD\eta\>=H.
\en

\bLemma\label{eta_nonzero}
We have
\eq
\rho\notin\cR \implies \eta \not = 0
\en
\eLemma
\bProof
See the proof of Proposition 2.1 in \cite{armony_questa03}.  In short, we have that
\eq
\rho\notin\cR \implies
\limsup_{t\go\infty}\frac{X_q(t)}{t} > 0 \mbox{ for some } q\in\cQ
\en
This in turn implies that:
\eq
\eta \not = 0 \mbox{ and, so }\limsup_{t\go\infty}\<\frac{X(t)}{t},\bD \frac{X(t)}{t}\> =
\<\eta,\Delta\eta\>=H > 0.
\en
\eProof

This result follows directly from the stability results of \cite{armony_questa03,ross_ton09} which say that when $\rho \in \cP$  the time-scaled  backlog goes to 0 as $t \go \infty$.  When the system is not stable, the opposite occurs and there exists a non-zero sub-limit. Our key result in this section is that $\eta$ is the \textbf{unique} limit of $\frac{X(t)}{t}$ over all possible arrival traces.

\bTheorem \label{th:limit}
When $\rho\notin\cR$, we have
\eq
\lim_{t\go\infty}\frac{X(t)}{t}= \eta\neq 0 .
\en
That is, the workload explodes on the same non-zero ray $\eta$ on any arrival trace.  Furthermore, $\eta$ is defined as the unique solution to the following convex program:
\eq\label{eq:cvxprog}
\<\eta,\bD\eta\> = \min_{\eta'\in\Psi(\rho,\cS)}\<\eta',\bD\eta'\>
\en
where
\eq
\Psi(\rho,\cS)=\{\eta': \eta'=(\rho - r)^+ \mbox{ with } r\in\cR\}
\en
and $\cR$ is the stability region given by $\cS$. Therefore, $r=\sum_{S\in\cS}\alpha_S S$ with $\sum_{S\in\cS}\alpha_S \le 1$ and $\alpha_S \ge 0$ for each $S\in\cS$, where $\cS$ is the set of service vectors.  Equivalently, $\eta$ is the unique fixed point which satisfies $\eta = \rho - \sum_{S\in\cS} \alpha_S S$ with $\sum_{S\in\cS} \alpha_S = 1, \alpha_S > 0$ and
\eq
\<\eta, \bD S_m\> \geq \<\eta,\bD S_k\> \forall m\textrm{ such that }\alpha_m > 0
\en
\eTheorem
\bProof
The proof will conclude at the end of Section \ref{sec:asym_dyn}.  We will first show that the limit is unique for each individual arrival trace.  We then extend our analysis to show the limit is the same for all arrival traces with identical time-average traffic load, $\rho$, which in turn, implies uniqueness over all arrival traces.  We begin with a number of definitions and structural properties of the defined elements.

From Section V of \cite{ross_ton09} on MaxWeight Scheduling cone geometry, recall that
$\cC_S=\{x\in\reals^Q:\<S,\bD x\>=\max_{S'\in\cS}\<S',\bD x\>\}$ is a cone, and when $X(t)\in\cC^o_S$ (the interior of $\cC_S$) the MaxWeight scheduling policy will choose $S(t)=S$. Moreover, the {\em surrounding cone} of any non-zero vector $\eta$ is the cone
\eq
\cC(\eta)=\bigcup_{S\in\cS^*(\eta)-\cS^\dag}C_S
\en
where $S^*(\eta)={\rm argmax}_{S\in\cS}\<S,\bD\eta\>$ is the set of service vectors of that MaxWeight Scheduling would select for backlog $\eta$ and $\cS^\dag$ is the set of {\em non-essential} ones (see \cite{ross_ton09}, end of Section IV).  Loosely speaking, the non-essential service vectors are the one whose removal will not change the stability region. Hence, $\cC(\eta)$ is the union of all cones corresponding to the set of service vectors which can be used if $X(t) = \eta$.  We have
\eq\label{tata}
X(t)\in\cC^o(\eta)\implies \<S(t),\bD \eta\> = \max_{S\in\cS}\< S,\bD\eta\> 
\en
where $\cC^o(\eta)$ is the interior of $\cC(\eta)$. Define now
\eq
\cK(\eta)=\{x\in\reals_{0+}^Q: x_q > \max_{S\in\cS}\{S_q\}
\mbox{ for each } q \mbox{ with } \eta_q > 0\},
\en
which is the set of backlogs such that no queue $q$ with $\eta_q > 0$ can be emptied in a single time slot. Note that $\cK(\eta)$ is upward-scalable; indeed, $x\in\cK(\eta)$ implies $\alpha x\in\cK(\eta)$ for any scalar $\alpha > 1$. Note that when $X(t)\in\cK(\eta)$ we have $X_q(t)>\max_{S\in\cS}\{S_q\}$ for all $q\in\cQ$ with $\eta_q>0$, so $D_q(t)=\min\{X_q(t),S_q(t)\}=S_q(t)$. Therefore,
\eq\label{zaza}
X(t)\in\cK(\eta)\implies D_q(t)=S_q(t) \mbox{ for all } q\in\cQ \mbox{ with } \eta_q > 0 .
\en
That is, all service capacity allocated at slot $t$ to queue $q$ with $\eta_q>0$ is used; there is no idling in that time slot; and the departures from all queues $q$ with $\eta_q >0$ is exactly equal to the total service provided to that queue. Consider now the set
\eq
\label{lala}
\cV(\eta)=\cK(\eta)\bigcap\cC^o(\eta)
\en
and note that it is upward-scalable, that is, $x\in\cV(\eta)$ implies $\alpha x \in\cV(\eta)$ for any scalar $\alpha > 1$. Thus, the set $\cV(\eta)$ is `cone-like'.

The main idea for the proof for the existence of our limit is that excursions away from $\eta$ take a very long time--so long that the time-average properties we have for arrivals will become active.  In particular, there will be some time after which all deviations of $\frac{X(t)}{t}$ away from $\eta$ will remain in $\cV(\eta)$.  To show this, we will need to demonstrate certain properties of these excursions.

\subsection{Structural Properties}\label{ssec:strucprop}
We now identify a number of structural properties of $\eta$ and related elements.  The proofs of these properties can be found in the Appendix.  These properties are important to characterizing deviations from $\eta$ and showing our main result on the limit of $\frac{X(t)}{t}$.

\bLemma\label{lm:str1} For every sequence $\{t'_c\}$ such that $t'_c<t_c$ and
\eq
X(t)\in\cV(\eta) \mbox{~~ for every ~~} X(t)\in (t'_c,t_c]
\en
for every $c$, we have
\eq
\<\frac{X(t_c)-X(t'_c)}{t_c-t'_c},\bD\eta\> =
\<\frac{\sum_{t=t'_c}^{t_c-1}A(t)}{t_c-t'_c},\bD\eta\> - \max_{S\in\cS}\<S(t),\bD\eta\>.
\en
\eLemma
During the interval $(t'_c,t_c]$, $X(t)\in \cK(\eta)$, which means any service vector used in this interval will maximize the inner product: $\<S(t),\bD\eta\>$.  Furthermore, by definition of $\cK(\eta)$, there is no idling for all $q$ such that $\eta_q > 0$.  This result simply accounts for the new jobs which arrive and the jobs which are serviced over the subset of queues with  $\eta_q > 0$.

\bLemma\label{lm:str2}
For any increasing unbounded time sequences $\{t_n\}$ and $\{t'_n\}$, we have
\eq
\lim_{n\go\infty}\frac{t_n-t'_n}{t_n}=\chi\in (0,1] \implies
\lim_{n\go\infty}\frac{\sum_{t=t'_n}^{t_n-1}A(t)}{t_n-t'_n} = \rho
\en
\eLemma
Therefore we can define the time-average workload over sub-intervals which grow linearly in time.

\bLemma\label{lm:str3} For any increasing unbounded subsequence $\{t_m\}$ with $\lim_{m\go\infty}\frac{X(t_m)}{t_m}=\mu$, we have
\eq
\<\mu,\Delta\eta\> \ge
\<\rho,\bD\eta\> - \max_{S\in\cS}\<S,\bD \eta\>
\en
\eLemma

\bLemma\label{lm:str4} For any increasing unbounded subsequence $\{t_m\}$ with $\lim_{m\go\infty}\frac{X(t_m)}{t_m}=\mu$, we have
\eq
\<\mu,\bD\eta\> \ge \<\eta,\bD\eta\> \implies\mu=\eta
\en
\eLemma

\bLemma\label{lm:str5}
For every $\epsilon\in (0,1)$ we have
\eq
- \left[\< \rho, \bD\eta \> - \max_{S\in\cS} \< S, \bD\eta \> \right]
\frac{\epsilon}{1-\epsilon} + \<\eta,\bD\eta\> \frac{1}{1-\epsilon}
\ge \<\eta,\bD\eta\>
\en
\eLemma

\subsection{Uniqueness of limit  $\mathbf{\lim_{t\go\infty}\frac{X(t)}{t}}$ on an individual arrival trace}\label{sec:uniquelimit}
In order to prove the uniqueness of the limit on a given arrival trace, we need one more definition and one more property.  We shall begin by assuming that there is some other convergent subsequence  $\{X(t_a)\}$ such that $\lim_{a\go\infty}\frac{X(t_a)}{t_a}=\psi\neq \eta$. Note that $\psi_q < \infty$ for all $q$.  This is easy to see since $\psi_q = \lim_{a\go\infty}\frac{X_q(t_a)}{t_a} \leq\lim_{a\go\infty}\frac{A_q(t_a)}{t_a} = \rho_q < \infty$.  We will eventually show that this subsequence does not exist.

Recall that we have subsequence $\{t_c\}$ which achieves the `$\limsup$' in \eqref{eq:def_eta}:
\eq
\lim_{c\go\infty}\frac{X(t_c)}{t_c} = \eta
\en
Define
\eq
s_c=\max\{t_a: t_a < t_c\} < t_c
\en
Since $s_c$ is a subsequence of $t_a$, $\lim_{c\go\infty}\frac{X(s_c)}{s_c} = \psi$.  This is a subsequence of the deviations away from $\eta$.
Given this definition of $\{s_c\}$, we can show the following property:
\bLemma\label{lm:eps_bounds}
We have that
\eq \label{eq:liminf}
\liminf_{c\go\infty}\frac{t_c-s_c}{t_c}=\epsilon\in(0,1)
\en
\eLemma
i.e. $s_c$ grows nearly linearly with $t_c$. The proof of this Lemma can be found in the Appendix.

We are now prepared to show the uniqueness of the limit on an individual arrival trace.  That is:
\bProposition\label{prop:uniquelimit}
There is no subsequence $\mathbf{\{t_a\}}$ with $\mathbf{\lim_{a\go\infty}\frac{X(t_a)}{t_a}=\psi\ne\eta}$.  Therefore,
\eq
\lim_{t\go\infty}\frac{X(t)}{t} = \eta
\en
\eProposition
\bProof
Arguing by contradiction, assume that there is some other convergent subsequence  $\{X(t_a)\}$ such that $\lim_{a\go\infty}\frac{X(t_a)}{t_a}=\psi\neq \eta$. We shall show that this is impossible. As stated before, note that $\psi_q < \infty$ for all $q$.
Select now a subsequence of $\{t_c\}$ on which the `$\liminf$' is attained in \eqref{eq:liminf},
but keep the same indexing $c$ of the original one for notational simplicity, hence,
\eq\label{1a1a}
\lim_{c\go\infty}\frac{t_c-s_c}{t_c}=\epsilon\in (0,1) .
\en
by Lemma \ref{lm:eps_bounds}; hence the length of time of the deviations  grows linearly with $t_c$. Therefore, $\lim_{c\go\infty}\frac{t_c}{s_c}=\frac{1}{1-\epsilon}$ and
$\lim_{c\go\infty}\frac{t_c-s_c}{s_c}=\frac{\epsilon}{1-\epsilon}$.

Applying Lemmas \ref{lm:str1} and \ref{lm:str2} with $\{t'_c\}=\{s_c\}$, dividing by $t_c-s_c$ and letting $c\go\infty$, we get
\eq\label{dada}
\lim_{c \rightarrow \infty}
\<\frac{X(t_c)-X(s_c)}{t_c-s_c}, \bD \eta \>
=  \< \rho, \bD\eta \> - \max_{S\in\cS} \< S, \bD\eta \> .
\en
Then, we can write
\begin{eqnarray}\label{key-ineq}
\< \psi, \bD \eta\> &=&
\lim_{c\go\infty}
\<\frac{X(s_c)}{s_c}, \bD \eta\>\nonumber \\
&=&
\lim_{c\go\infty}
\< - \frac{X(t_c)-X(s_c)}{t_c - s_c} \ \frac{t_c-s_c}{s_c}
+ \frac{X(t_c)}{t_c}\ \frac{t_c}{s_c} ,\bD \eta \> \nonumber \\
&=& - \left[\< \rho, \bD\eta \> - \max_{S\in\cS} \< S, \bD\eta \> \right]
\frac{\epsilon}{1-\epsilon} + \<\eta,\bD\eta\> \frac{1}{1-\epsilon}\nonumber \\
&\ge& \<\eta,\bD\eta\>
\end{eqnarray}
The last equality is due to Lemma \ref{lm:str5}. Therefore,
$\<\psi,\bD\eta\> \ge \<\eta,\bD\eta\>$, which implies $\psi=\eta$ from Lemma \ref{lm:str4}. But this contradicts the assumption that $\psi\ne\eta$. This establishes the sought after contradiction. So for each individual arrival trace, there exists a unique limit $\lim_{t\go\infty}\frac{X(t)}{t} = \eta$, which concludes the proof of the proposition. Moreover, since $\lim_{t\go\infty}\frac{X(t)}{t} = \eta$ this implies that there exists $t_o < \infty$ such that $X(t)$ is in $\cV(\eta)$ for all $t>t_o$.
\eProof

The main argument is essentially that no other sublimit can get very far from $t_c$ and any excursion away from $\eta$ is small.  Therefore at large $t$, all subsequences are close enough to $\eta$ to be in one of the neighboring cones, i.e. in $\cV(\eta)$.  We have now shown that there exists a unique limit, $\eta$, such that on a given arrival trace: $\lim_{t\go\infty} \frac{X(t)}{t} = \eta$.  To complete the proof of Theorem \ref{th:limit}, it remains to show that the limit, $\eta$, is independent of the particular arrival trace.  To do this, we turn to characterizing $\eta$.

\subsection{Characterizing the limit $\eta$}\label{ssec:eta}
The purpose of this section is to characterize the limit $\eta$ in terms of $\rho$ and service vectors $\cS$ to establish the independence of $\eta$ on the individual arrival trace.   Knowing that $\lim_{t\go\infty}\frac{X(t)}{t}=\eta$, we now turn to identifying a couple of the characteristic properties of $\eta$.  The proofs of these Lemmas can be found in the Appendix.

\bLemma\label{fixedpoint}
Every limit is a fixed point.  That is,
\eq
\eta = \lim_{t\rightarrow \infty} \frac{X(t)}{t} = \Big [ \rho - \sum_{m=1}^N \alpha_m S_m \Big ]^+
\en
for some $\alpha_m \geq 0, \sum_m \alpha_m = 1$.  Furthermore, $\alpha_m > 0$ implies that $\eta \in C_{S_m}$.
\eLemma
Note that because $\alpha_m>0$ implies that $\eta \in C_{S_m}$, we know that $\eta$ is in the intersection of all cones with $\alpha_m>0$.  Hence, if there are multiple $\alpha_m > 0$, then $\eta$ is on the boarder of all  cones with $\alpha_m > 0$.  Because $\eta$ is a fixed point:
\bLemma \label{lemma:etarho}
The following equality holds:
\eq
\<\eta,\bD\eta\> = \<\rho,\bD\eta\> - \max_{S\in\cS}\<S,\bD\eta\>
\en
\eLemma

We are now in position to show that $\eta$ depends only on the load vector, $\rho$, and available service vectors, $\cS$.  That is, $\eta$ is independent of the particular arrival trace and is the solution of a simple convex program.
\bProposition\label{prop:etaunique}
The vector $\eta=\lim_{t\go\infty}\frac{X(t)}{t}$ is the \underline{unique} minimizer of
\eq\label{eq:cprog}
\<\eta,\bD\eta\> = \min_{\eta'\in\Psi(\rho,\cS)}\<\eta',\bD\eta'\>
\en
where
\eq
\Psi(\rho,\cS)=\{\eta': \eta'=(\rho - r)^+ \mbox{ with } r\in\cR\}
\en
and $\cR$ is the stability region given by $\cS$. Therefore, $r=\sum_{S\in\cS}\alpha_S S$ with $\sum_{S\in\cS}\alpha_S \le 1$ and $\alpha_S \ge 0$ for each $S\in\cS$, where $\cS$ is the set of service vectors.
\eProposition
\bProof
First we show that $\eta$ is indeed the solution to \eqref{eq:cprog}.  Then we show there is only one solution in order to conclude that $\eta$ is unique.

From Lemma \ref{fixedpoint}, we have that $\eta \in \Psi(\rho,\cS)$.  Arbitrarily choose any vector
\eq
\bar{\eta}=\Big (\rho - \sum_{S\in\cS}\alpha_S S\Big )^+, \mbox{  with  }
\sum_{S\in\cS}\alpha_S \le 1,  \mbox{  and  } \alpha_S\ge 0, S\in \cS .
\en
Projecting on $\bD\eta$ we get
\begin{eqnarray}
\<\bar{\eta},\bD\eta\> &=&\< \Big [\rho - \sum_{S\in\cS}\alpha_S S\Big ]^+,\bD\eta\> \nonumber \\
 &\geq& \< \rho - \sum_{S\in\cS}\alpha_S S ,\bD\eta\> \nonumber \\
 & = &\<\rho,\bD\eta\> - \sum_{S\in\cS}\alpha_S \<S,\bD\eta\> \nonumber \\
&\ge& \<\rho,\bD\eta\> -\max_{S\in\cS}\<S,\bD\eta\> = \<\eta,\bD\eta\>
\end{eqnarray}
The first inequality comes from the fact that $\bD_{qq} >0$ and $\eta_q \geq 0$. The last equality comes from Lemma \ref{lemma:etarho}.  Therefore, $\<\bar{\eta},\bD\eta\> \ge \<\eta,\bD\eta\>$. This implies (recalling that $\bD$ is positive-definite) that
\begin{eqnarray}
0&\le& \<\bar{\eta}-\eta,\bD(\bar{\eta}-\eta)\> \nonumber \\
&=&
\<\bar{\eta},\bD\bar{\eta}\> -2\<\bar{\eta},\bD\eta\>+\<\eta,\bD\eta\> \nonumber \\
&\le& \<\bar{\eta},\bD\bar{\eta}\> -2\<\eta,\bD\eta\>+\<\eta,\bD\eta\> \nonumber \\
&=&
\<\bar{\eta},\bD\bar{\eta}\> - \<\eta,\bD\eta\>,
\end{eqnarray}
so $\<\bar{\eta},\bD\bar{\eta}\> \ge \<\eta,\bD\eta\>$ and $\eta$ is the minimizer of $\<\eta',\bD\eta'\>$.

We still need to prove that the minimizer $\eta$ is unique. This is done by showing that 1) $\<\eta,\bD\eta\>$ is strictly convex in $\eta$ and 2) the set $\Psi(\rho,\cS)$ is convex--uniqueness will follow from convex programming theory.  1) It is trivial to show that $\<\eta,\bD\eta\>$ is strictly convex in $\eta$ since $\bD > 0$ is a positive definite matrix.  2) We now show that the set $ \Psi \equiv \Psi(\rho,\cS) $ is convex.  First, we see that for any $r\in\cP$ and corresponding $x = (\rho - r)^+ \in \Psi$, there exists $\bar{x} = \rho - \bar{r} = (\rho - r)^+ = x$ with $\bar{r} \in \cP$.  Let $\bar{r}_k = \min(\rho_k,r_k) \leq r_k$.  Since $\bar{r} \leq r \in \cP$, then $\bar{x} \in \Psi$.  Now, consider two vectors $x,x' \in \Psi$ with corresponding $r,r' \in \cP$ such that $x = (\rho - r)^+$ and $x' = (\rho - r')^+$.  What remains to be shown is that for any $a \in [0,1]$, $ax+(1-a)x' \in \Psi$.  Indeed, we have:
\begin{eqnarray}
ax+(1-a)x' &= &a\bar{x} + (1-a)\bar{x}' \nonumber \\
& = & a(\rho - \bar{r}) + (1-a)(\rho -\bar{r}') \nonumber \\
& = & \rho - (a\bar{r} + (1-a)\bar{r}')
\end{eqnarray}
By the convexity of $\cP$, we know that $a\bar{r} + (1-a)\bar{r}' \in \cP$ and subsequently, $ \rho - (a\bar{r} + (1-a)\bar{r}') \in \Psi$.  This concludes the proof.
\eProof

Based on the characterization of $\eta$ as the solution to the convex program \eqref{eq:cprog}, via convex optimization theory, we can conclude that there is only one fixed point.
\bLemma \label{lemma:onefixedpt}
There exists exactly one fixed point, $\eta$:
\begin{eqnarray}
\eta &=& [\rho - \sum_m\alpha_m S_m]^+ \nonumber \\
 0&\leq&\alpha_m  \nonumber \\
1 &=& \sum_m \alpha_m\nonumber \\
\alpha_m > 0 &\implies &\<\eta, \bD S_m\> \geq \< \eta, \bD S_k\>, \forall k
\end{eqnarray}
\eLemma

We have just shown that on all arrival traces with system load $\rho$, $\lim_{t\go\infty}\frac{X(t)}{t} =\eta$, is unique.  Furthermore, the limit, $\eta$, is identical across all such traces.  $\eta$ can be characterized as the unique solution to the convex program \eqref{eq:cvxprog}; equivalently, it is the unique fixed point of our system.  This concludes the proof of Theorem \ref{th:limit}.

To summarize, whenever the system is in overload the workload grows along a vector defined by $\eta$ which is the solution of the convex program in Proposition \ref{prop:etaunique}.  Moreover, $\eta$ is
{\em independent of the particular arrival trace}.  From Lemma \ref{fixedpoint}, we know that $\eta$ is on the intersection of some set of cones.  If there exists only one $\alpha_m > 0$, then $\eta = [\rho - S_m]^+$ and is in $C_{S_m}$.  If there are multiple $\alpha_m > 0$, then $\eta = [\rho - \sum_mS_m]^+$ is on the boundary of the  set of cones with $\alpha_m >0$.  We will utilize this information to control for our fairness criterion.

\section{Fair Control via the $\bD$ matrix}\label{sec:control}

We have now seen how the asymptotic behavior of the workload vector, $X(t)$, behaves given service vectors $\cS$ and MaxWeight matrix $\bD$.  In particular, when the system load is outside of the stability region, the workload will explode along a single vector.  During a long period of temporary stress, the queueing system  is effectively \emph{unstable} during this window and the valuable service resources become strained.  Under the MaxWeight scheduling policy, queues with exceptionally high load will starve resources from other, less stressed, queues.  This begs the question of how to share resources in a fair manner when the system is unstable.  Our goal in this section is to consider how to manipulate the $\bD$ matrix in order to ensure `fairness' in this queueing system.

As we have discussed before, there are many different definitions of fairness.  We focus on a notion of fairness where the workloads grow according to a fixed proportion.  More formally, we assume we are given a set of ratios $\theta_q\geq 0, \sum_q \theta_q = 1$.  Whenever $\rho \not \in \cP$,  the goal is to control the workload such that:
\eq
\lim_{t\rightarrow \infty}\frac{X_q(t)}{\sum_k X_k(t)} = \theta_q, \forall q \in \cQ
\en
From Theorem \ref{th:limit}, we know that $\lim_{t \rightarrow \infty}\frac{X(t)}{t} = \eta$ and we want:
\eq
\frac{\eta_q}{\sum_i \eta_i} = \theta_q, \forall q
\en
When $\rho \in \cP$, the system is stablizable, and so the workloads should remain finite. For large $t$:
\eq
\frac{X_q(t)}{t} \rightarrow 0, \forall q
\en
We know that the MaxWeight scheduling policy guarantees the second criterion--it is a stabilizing policy \citep{ross_ton09}.  The focus of this section is to show that the MaxWeight scheduling policy also satisfies the first criterion under certain conditions via proper specification of the $\bD$ matrix.

\subsection{Control via the $\bD$ matrix}\label{ssec:Bcontrol}
From Lemma \ref{fixedpoint}, we know that $\eta$ is on the intersection of the set of cones with $\alpha_m > 0$ in the definition of $\eta = (\rho-\sum_m \alpha_m S_m)^+$.  This corresponds is a cone boundary if there are more than one $\alpha_m > 0$. This boundary depends on the $\bD$ matrix used in the MaxWeight scheduling policy.   We now consider how we can choose the $\bD$ matrix to place cone boundaries and, subsequently, achieve the desired fairness criterion.

\begin{lemma}\label{lemma:bplace} In $Q$-dimensions, consider any $M$ service vectors $S_{i_1}, S_{i_2}, \dots, S_{i_M}$.  Suppose there exists a diagonal positive definite matrix $\hat{\bD}$ and non-negative $Q$-dimensional vector ${v}\geq0$ such that ${v}$ is a \emph{boundary vector} of the $M$ cones, i.e. for each $k \in [1,M]$ and for all $j$:
\eq\label{eq:boundary}
\<{v},\hat{\bD}S_{i_k}\> \geq \<{v},\hat{\bD}S_j\>
\en
Then, for any non-negative vector $\eta$ such that $\eta_q = 0$ if and only if $v_q = 0$, there exists a diagonal positive definite matrix $\bD$ such that for each $k \in [1,M]$ and all $j$:
\eq
\<\eta,\bD S_{i_k}\> \geq \<\eta,\bD S_j\>
\en
i.e. the boundary between the $M$ cones can be placed arbitrarily in $\reals^Q_+$.  This matrix is specified as:
\eq
\bD_{qq} = \left\{
             \begin{array}{ll}
               \frac{\hat{\bD}_{qq}{v}_q}{\eta_q}, & \hbox{$v_q >0$;} \\
               1, & \hbox{$v_q = 0$.}
             \end{array}
           \right.
\en
\end{lemma}
This proof is given in the Appendix.

Because Lemma \ref{lemma:bplace} holds for \emph{any} positive definite diagonal matrix $\hat{\bD}$, it will hold for $\hat{\bD} =  \bI$.  Furthermore, it is easy to verify condition \eqref{eq:boundary} when $\hat{\bD} = \bI$, as this results in taking simple inner products.
\bCorollary\label{cor:iplace}
The boundary, $\eta \geq 0$, between $M$ cones defined by service vectors $S_{i_1}, S_{i_2}, \dots, S_{i_M}$ can be placed arbitrarily in $\reals^Q_+$ by using MaxWeight scheduling with matrix $\bD$ if there exists $v \geq 0$ on the boundary of the $M$ cones induced by $\hat{\bD} = \bf{I}$:
\eq
\<v,S_{i_k}\> \geq \<v, S_j\>, \forall k\in[1,M], \forall j
\en
Furthermore,  $v_q = 0$ if and only if $\eta_q = 0$. The required $\bD$ is:
\eq
\bD_{qq} = \left\{
             \begin{array}{ll}
               \frac{{v}_q}{\eta_q}, & \hbox{$v_q >0$;} \\
               1, & \hbox{$v_q = 0$.}
             \end{array}
           \right.
\en
\eCorollary

Under certain necessary and sufficient conditions we can arbitrarily place the cone boundary and we can then control the workload to explode along a desired vector defined by $\eta$.
\bProposition\label{prop:control}
There exists MaxWeight matrix $\bD$ such that
\eq
\lim_{t\go\infty}\frac{X(t)}{t} = \eta
\en
if and only if the following conditions hold:
\begin{enumerate}
  \item $\eta = (\rho - \sum_m\alpha_m S_m)^+$ for some $\alpha_m \geq 0, \sum_m \alpha_m = 1$.
  \item There exists $v \geq 0$ such that $\alpha_m > 0$ implies $\<v, S_m\> \geq \<v, S_k\>$ for all $k$.  Furthermore,  $v_q = 0$ if and only if $\eta_q = 0$.
\end{enumerate}
\eProposition

This allows us to characterize the {\em feasible} fairness criterion:
\bCorollary\label{cor:feasible}
Fairness criterion $\theta$ is feasible via MaxWeight Matching if and only if:
 \begin{enumerate}
   \item There exist $M$ service vectors and $v\geq0$, with $v_q = 0$ if $\theta_q = 0$, such that
\eq
\<v,S_{i_k}\> \geq \<v, S_j\>,\forall k\in[1,M], \forall j
\en
\item For vectors $S_{i_m}$ in 1, there exist $\alpha_m \geq 0$, $\sum_{m=1}^M \alpha_m = 1$ such that $\eta = (\rho - \sum_{m=1}^M \alpha_mS_{i_m})^+$.  Additionally, for each $\alpha_m >0$:
\eq
\<v,S_m\> \geq \<v,S_k\>,\forall k
\en
\item $\eta$ must satisfy:
\eq
\theta_q = \frac{\eta_q}{\sum_k \eta_k}
\en
 \end{enumerate}
 \eCorollary
This is a direct consequence of Proposition \ref{prop:control}.

\subsubsection{Feasible Criteria: An Example} 
\begin{figure}[ht]
\centering
\subfigure[Stability Region]{
\includegraphics[scale=.4]{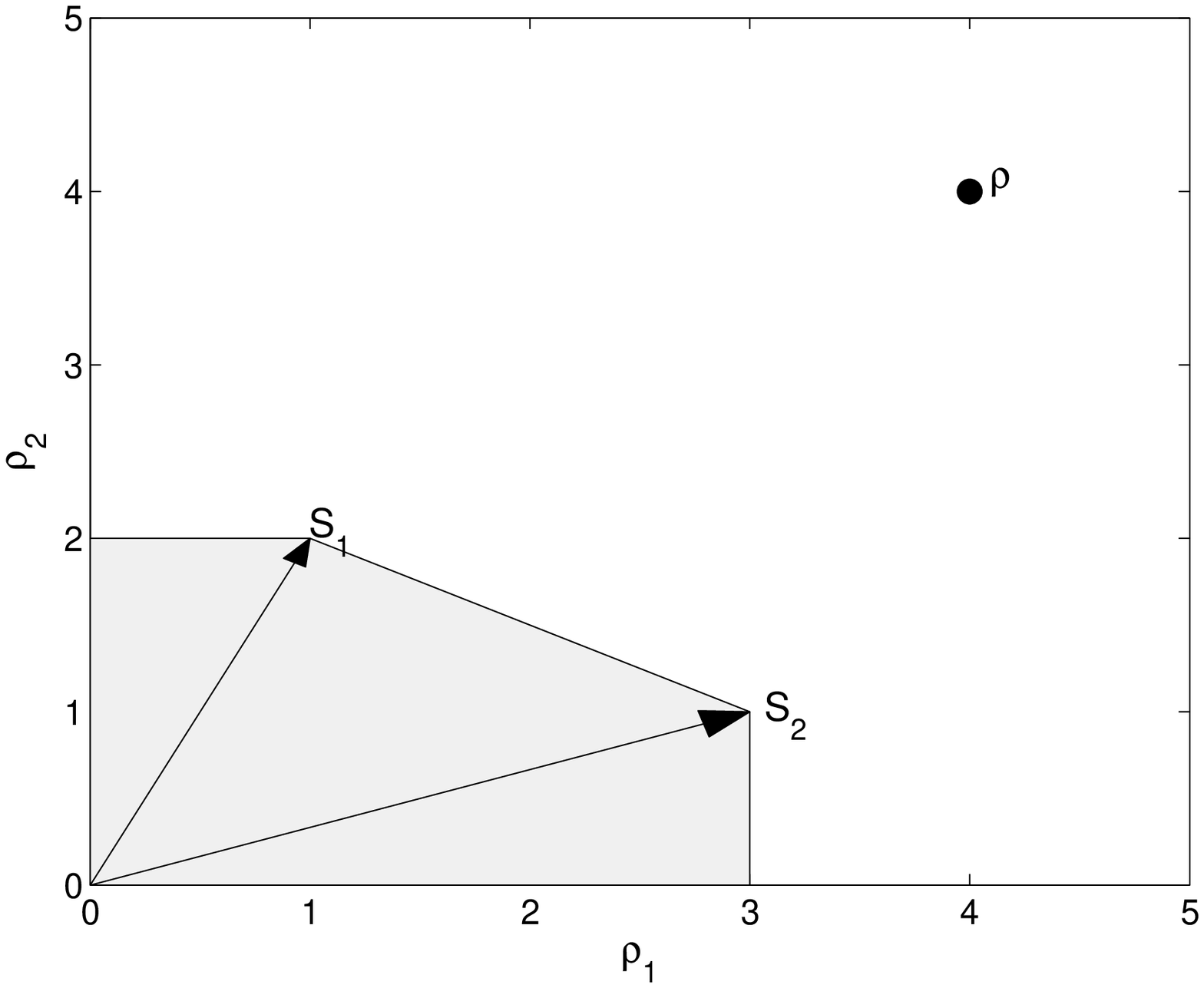}
}
\subfigure[Feasible Fairness Criteria]{
\psfrag{th}{$\hat{\theta}$}
\psfrag{thh}{$\theta$}
\includegraphics[scale=.4]{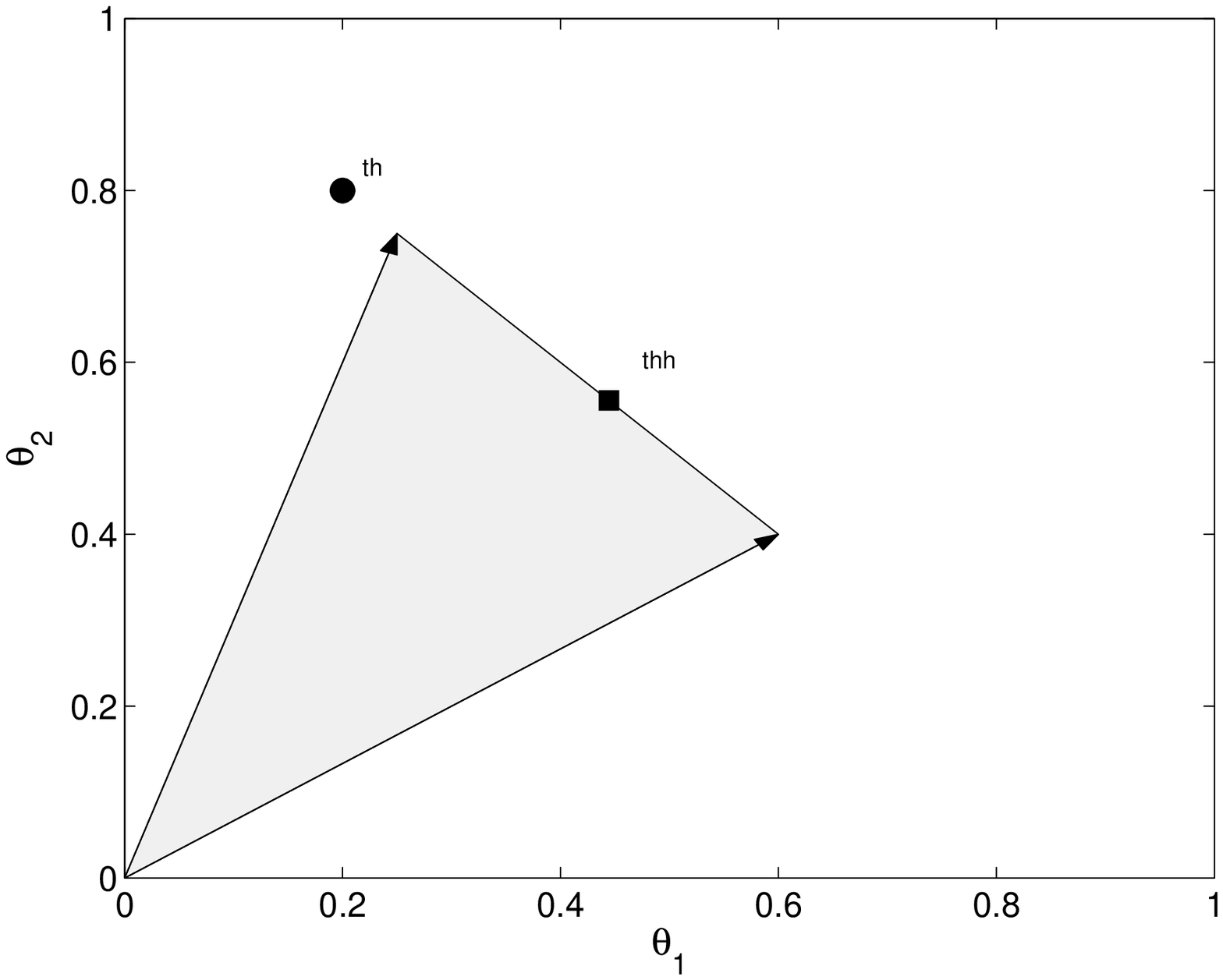}
}
\label{fig:feas}
\caption{Feasible Fairness Criteria for $\rho$: $N = 2$ service vectors and $Q = 2$ queues.}
\end{figure}

We now present an example with  $N = 2$ service vectors and $Q = 2$ queues and specify the feasible criteria given $\rho = [4,4]$, $S_1 = [1,2]$ and $S_2 = [3,1]$.  It is easy to see that $v = [1,2]$ satisfies Condition 1 in Corollary \ref{cor:feasible} as long as $\theta > 0$.  Now, by Condition 2,
\begin{eqnarray}
\eta & = & (\rho - \alpha S_1 - (1-\alpha)S_2)^+ \nonumber \\
& = & \alpha(\rho - S_1) + (1-\alpha)(\rho - S_2) \nonumber \\
& = & \alpha\left[
  \begin{array}{c}
    1 \\
    3 \\
  \end{array}
\right] + (1-\alpha)\left[
  \begin{array}{c}
    3 \\
    2 \\
  \end{array}
\right], \forall \alpha \in [0,1]
\end{eqnarray}
To find the feasible $\theta$, we normalize $\eta$ as in Condition 3.  Hence, any $\theta = \alpha\left[
  \begin{array}{c}
    1/4 \\
    3/4 \\
  \end{array}
\right] + (1-\alpha)\left[
  \begin{array}{c}
    3/5 \\
    2/5 \\
  \end{array}
\right], \alpha \in [0,1]$ is feasible.   In Figure \ref{fig:feas}, we can see the stability region and $\rho$, which is outside of the stability region.  All fairness directions in the gray portion, such as $\theta$, are feasible.  Other fairness criteria, such as $\hat{\theta}$, are infeasible.

\subsubsection{An Infeasible Feasible Criterion} \label{ssec:infeasfeas}
One might naturally ask whether the set of feasible fairness directions according to MaxWeight Matching includes all feasible directions when a more general set of policies is allowed. It turns out the  answer is no. We demonstrate this via an example where the fairness criterion is \emph{not} achievable via MaxWeight Matching but there exists a policy which does achieve it. Consider a scenario with $Q = 3$ queues and $N = 3$ service vectors:
\eq
S_1 = \left(
  \begin{array}{c}
    1 \\
    0 \\
    1 \\
  \end{array}
\right),
S_2 = \left(
  \begin{array}{c}
    0 \\
    1 \\
    1 \\
  \end{array}
\right),
S_3 = \left(
  \begin{array}{c}
    3/4 \\
    3/4 \\
    2 \\
  \end{array}
\right)
\en
Let the fairness criteria and system load be:
\eq
\theta = \left(
  \begin{array}{c}
    1/3 \\
    1/3 \\
    1/3 \\
  \end{array}
\right),
\rho = \left(
  \begin{array}{c}
    13/8 \\
    13/8 \\
    5/2 \\
  \end{array}
\right)
\en
It is easy to see that using a combination of all three service vectors, the fairness criteria can be achieved:
\eq
\eta = \left[\rho - \frac{1}{4}S_1 - \frac{1}{4}S_2 - \frac{1}{2}S_3\right]^+ = \left(
  \begin{array}{c}
    1 \\
    1 \\
    1 \\
  \end{array}
\right)
\en
However, this fairness criteria \emph{cannot} be met using MaxWeight Matching.  In particular, by Conditions 1 and 2 of Corollary \ref{cor:feasible}, there must exists some $v>0$ (since $\theta > 0$) such that:
\eq \label{eq:boarder}
\<v,S_1\> = \<v,S_2\> = \<v,S_3\>
\en
With some algebra, one can see that the $v$ which satisfies \eqref{eq:boarder} is
\eq
v = \gamma\left(
  \begin{array}{c}
    2 \\
    2 \\
    -1 \\
  \end{array}
\right)
\en
for any $\gamma$.  Hence, there does not exists $v >0$ and Conditions 1 and 2 of Corollary \ref{cor:feasible} cannot be satisfied.  This $\theta$ is not feasible via MaxWeight Matching.

We find that MaxWeight Matching allows us to control for a general, characterizable set of fairness criterion; however, there may be other policies which can achieve other criteria.  Note that as long as MaxWeight policies can achieve the fairness criteria, it is guaranteed to achieve it in the most efficient manner by Theorem \ref{th:perf}.

\subsection{Workload Minimization}\label{ssec:MWMfair}
Recall that our ultimate goal is to minimize the long run average backlog \eqref{eq:min} subject to the fairness criterion \eqref{eq:faircrit}.  So far we have established that, given a feasible fairness criterion $\theta$, the MaxWeight scheduling policy will satisfy this criterion with the right choice of the matrix $\bD$.  In particular, we have identified the  `direction' which the workload will grow as well as a methodology to control this direction.  However, we have  not specified the rate at which the workload will grow.  We have used the MaxWeight scheduling policy to achieve our fairness criterion and justified the use of these algorithms because they are stabilizing when $\rho \in \cP$, they are simple to implement, and they allow us to achieve our fairness criterion.  We now discuss another feature of the MaxWeight scheduling policy which makes it highly desirable for our purposes.  Namely, it achieves the \emph{smallest} workload which satisfies our fairness criterion.

\bTheorem\label{th:perf}
Let $X(t)$ be the workload achieved under the MaxWeight scheduling policy.  Let $\bar{X}(t)$ be the workload achieved under some other algorithm  such that:
\eq
\lim_{t\rightarrow \infty} \frac{X(t)}{t} = \eta, \lim_{t\rightarrow \infty} \frac{\bar{X}(t)}{t} = \bar{\eta}
\en
where both achieve the desired fairness criterion:
\eq
\frac{\eta_q}{\sum_k \eta_k} = \frac{\bar{\eta}_q}{\sum_k \bar{\eta}_k} = \theta_q, \forall q
\en
Then, $\eta$ is the \emph{minimal} workload vector which can achieve the fairness criterion, $\theta$:
\eq
\eta \leq \bar{\eta}
\en
\eTheorem
The proof of this Theorem is given in the Appendix.

While there may be many algorithms which achieve our fairness criterion, $\theta$, MaxWeight scheduling is efficient in the sense that no other algorithm has smaller backlogs. Note that our notion of efficiency is different than that from \cite{Shah_Wischik}.  In their work, they show that MaxWeight is \emph{not} efficient.  However, they define efficiency as the total throughput of the system, i.e. they show that MaxWeight does not minimize $\sum_q X_q(t)$ over \emph{all} possible directions. This is quite different from what we show, which is that, given a direction $\theta$, then $\sum_q X_q(t)$ is minimized.

\subsection{Robustness with Respect to $\rho$}\label{ssec:robust}
Thus far, we have assumed that the load vector $\rho$ is known.  Under this assumption, we are able to select the necessary MaxWeight matrix, $\bD$, to achieve the desired fairness criterion, $\theta$, as long as it is feasible. Now we suppose that $\rho$ is unknown  and examine whether we are still able to choose $\bD$ to achieve the desired fairness criterion.  Throughout this discussion, we will assume that $N > 1$; otherwise there is no control and $\eta = (\rho - S)^+$ for all  $\rho$, irrespective of $\bD$.

\begin{figure}[ht]
\centering
\subfigure[Stability Region]{
\includegraphics[scale=.4]{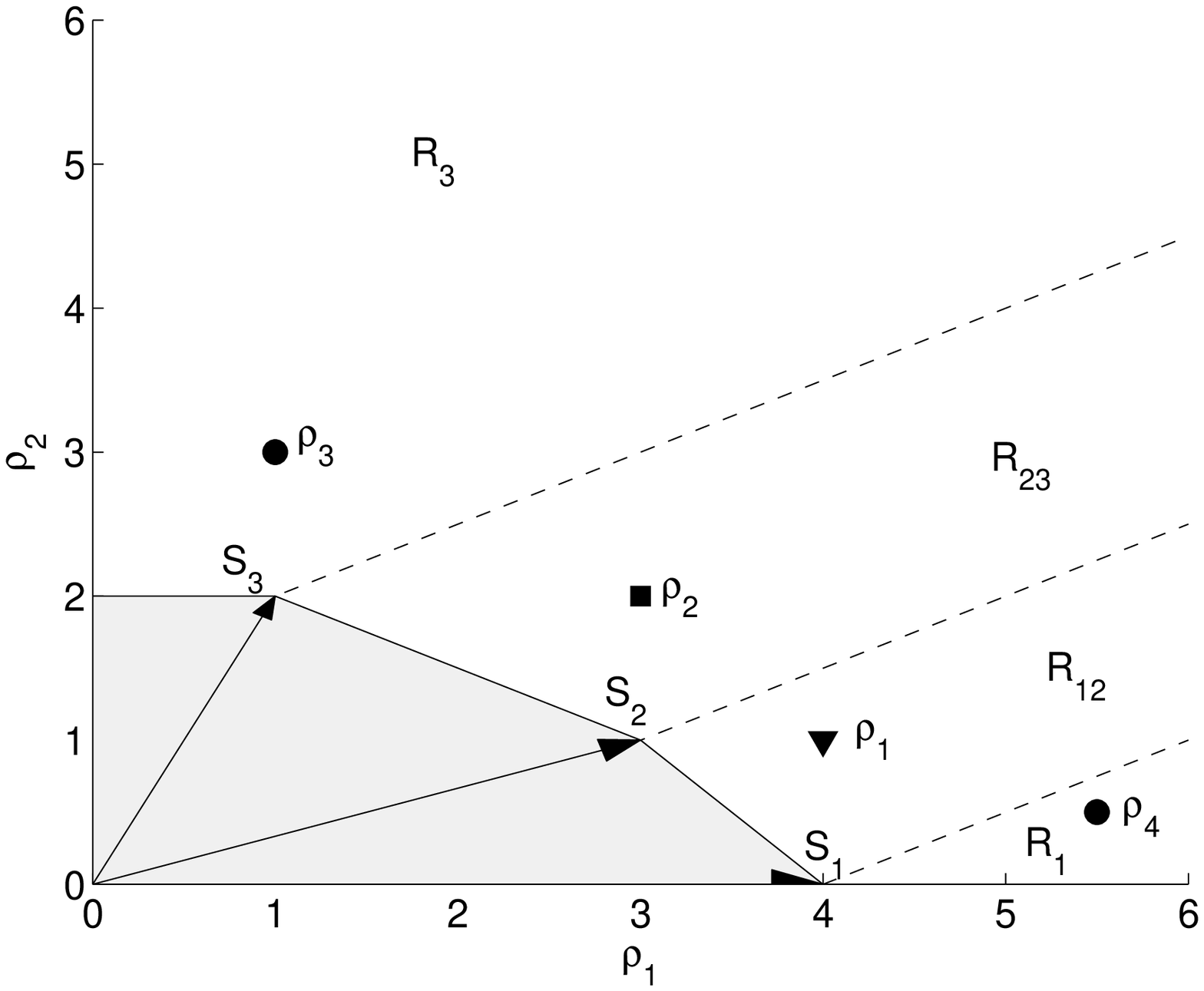}
\label{fig:st}
}
\subfigure[Scheduling cones]{
\includegraphics[scale=.4]{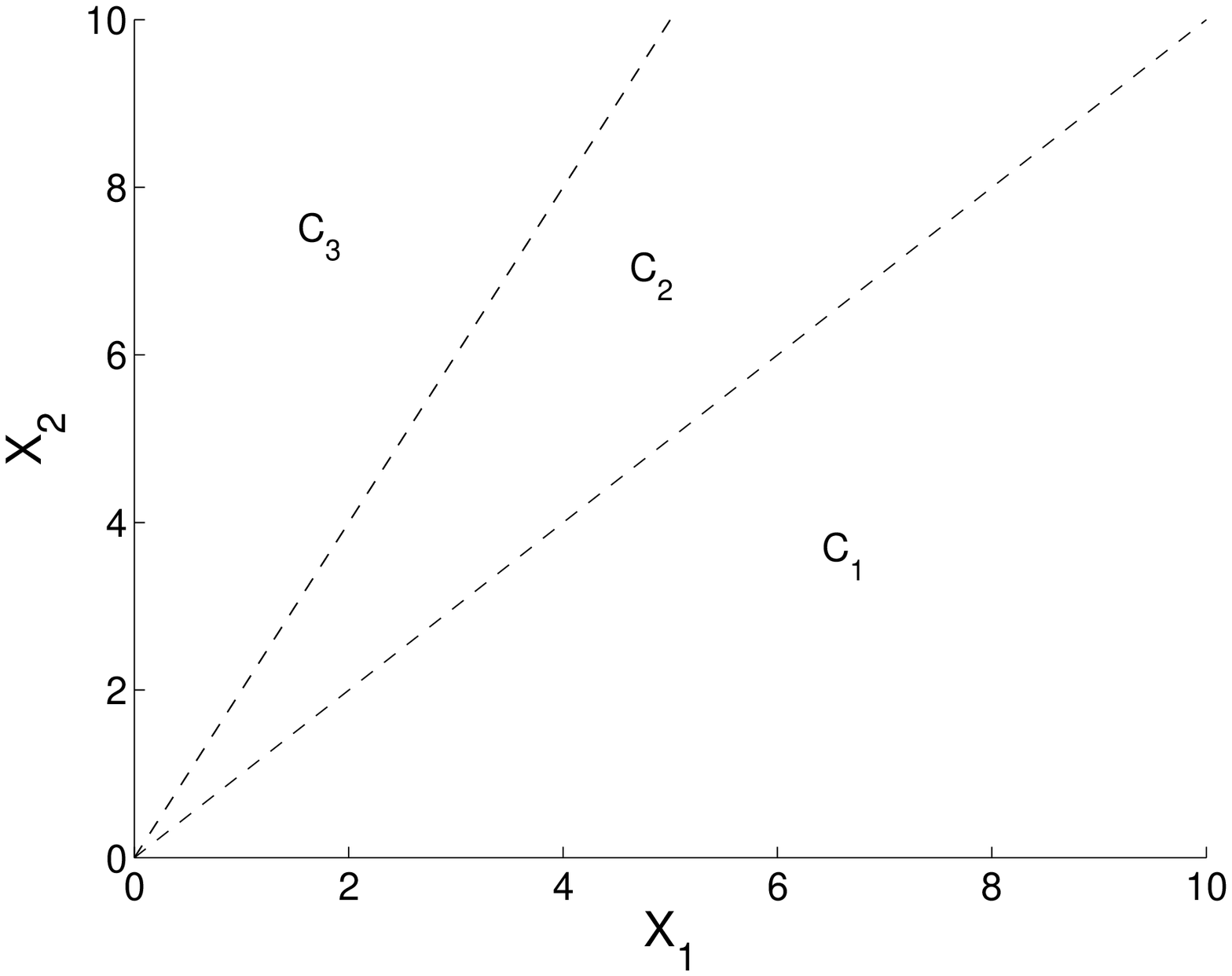}
\label{fig:sc}
}
\caption{Stability and Cone regions for $N = 3$ service vectors and $Q = 2$ queues.}
\label{fig:N3}
\end{figure}

Consider the following example with $N = 3$ service vectors and  $Q=2$ queues as depicted in Figure \ref{fig:N3}. Let
\eq
S_1 = \left(
  \begin{array}{c}
    4 \\
    0 \\
  \end{array}
\right),
S_2 = \left(
  \begin{array}{c}
    3 \\
    1 \\
  \end{array}
\right),
S_3 = \left(
  \begin{array}{c}
    1 \\
    2 \\
  \end{array}
\right)
\en

Suppose our fairness criteria is:
\eq
\theta = \left(
  \begin{array}{c}
    2/3 \\
    1/3 \\
  \end{array}
\right)
\en
We consider 3 different load vectors, which are outside of the stability region:
\eq
\rho_1 = \left(
  \begin{array}{c}
    4 \\
    1 \\
  \end{array}
\right),
\rho_2 = \left(
  \begin{array}{c}
    3 \\
    2 \\
  \end{array}
\right),
\rho_3 = \left(
  \begin{array}{c}
    1 \\
    3 \\
  \end{array}
\right),
\rho_4 = \left(
  \begin{array}{c}
    5 \\
    .5 \\
  \end{array}
\right)
\en
When the system load $\rho = \rho_3$ or $\rho_4$,Conditions 2 and 3 in Corollary \ref{cor:feasible} cannot be satisfied; hence, the fairness direction, $\theta$, is infeasible and there does not exist a MaxWeight matrix to achieve it.  With some algebra, we can see that the necessary MaxWeight matrix to achieve the fairness criteria depends on $\rho$:
\eq
\bD_1 = \left(
  \begin{array}{cc}
    1& 0 \\
    0 &2\\
  \end{array}
\right),
\bD_2 = \left(
  \begin{array}{cc}
    1& 0 \\
    0 &4\\
  \end{array}
\right)
\en
From Theorem \ref{th:limit},  the workload is given by $\eta=[\rho - \sum_m\alpha_mS_m]^+$. Hence, to achieve the desired fairness criteria, the goal is to find a point on the boundary of the stability region such that subtracting that point from the system load, $\rho$, results in a vector which is consistent with the fairness direction $\theta$.  From Figure \ref{fig:st}, we see that a point on the stability boundary which is given by the convex combination of $S_1$ and $S_2$ satisfies this constraint for $\rho_1$. The $\bD$ matrix skews the dimensions such that the boundary of interest is placed in the direction of the desired fairness criterion and makes this distance `minimal' as in Proposition \ref{prop:etaunique}.  Hence, the boundary vector of interest for $\rho_1$ is the boundary between cones $1$ and $2$.  This boundary vector can be moved to the fairness direction $\theta$ by using $\bD_1$.  Similarly, $\bD_2$ moves the boundary vector between cones $2$ and $3$ for $\rho_2$.  This example shows that the boundary vector of interest and, subsequently, the necessary MaxWeight matrix $\bD$ depends on $\rho$.

Despite the preceding example, it is possible to select $\bD$ without precise knowledge of $\rho$.  This ability depends on $k$, the number of subsets of service vectors $\{S_i\}$ of size greater than $1$ which satisfy:
\eq \label{eq:relevant}
\<v,S_{i_j}\>  = \<v,S_{i_n}\> > \<v,S_m\>, \forall S_{i_j},S_{i_n} \in \{S_i\}\textrm{ and } S_m \not \in \{S_i\}
\en
for some  $v \geq 0$, $v \not = \bf{0}$.  Hence, $v$ is a \emph{boundary vector} as it is a vector on the boundary between neighboring cones $\{S_i\}$. We refer to this boundary as a \emph{relevant boundary}.  Note that this boundary is potentially a hyperplane of dimension greater than $1$, so $v$ may not be unique.  Since $N >1$, there exists at least one  boundary which satisfies \eqref{eq:relevant} ($k \geq 1$).  Moreover,  with $N$ service vectors, there are $2^N - (N+1) $ subsets of size greater than 1.  Clearly, some subsets may not satisfy \eqref{eq:relevant}, so $k < 2^N -(N+1)$.  Our two cases are:

\begin{enumerate}
  \item $\bf{[k = 1]}$  In this case, there is exactly one subset of service vectors with a relevant boundary.  If the boundary of this subset is a hyperplane, then one can selected any boundary vector $v \geq 0$ (with $v_q = 0$ if and only if $\theta_q = 0$).  Because there is only one boundary, we can specify the MaxWeight Matrix, $\bD$, for any fairness criterion, $\theta$.  As long as the system load $\rho$ is such that $\theta$ satisfies Conditions 2 and 3 in Corollary \ref{cor:feasible}, this $\bD$ matrix will control the backlogs to grow along $\theta$.  Therefore, we can choose $\bD$ without knowing $\rho$.   $\rho$ simply determines whether a fairness direction $\theta$ is achievable.

  \item $\bf{[k>1]}$ In this case, there are  multiple (but finite) subsets of service vectors which satisfy \eqref{eq:relevant}. For each subset, $\cS(y) = \{S_i\}$ ($y \in [1,k]$), let $v(y)$ be a  non-negative boundary vector with $v(y)_q = 0$ if an only if $\theta_q = 0$.  There is a $\bD(y)$ such that $\bD(y)_{qq} = v(y)_q/\theta_q$ will place the boundary of the $\cS(y)$ cones at the direction specified by $\theta$.  If $\theta$ is feasible, there will be exactly one subset of service vectors $\cS(y)$ for each $\rho$ which satisfies Condition 2 of Corollary \ref{cor:feasible}.  There will be a range of $\rho$ corresponding to each subset $\cS(y)$, and subsequently $\bD(y)$. As long as our estimate of $\rho$ is within the accuracy of these ranges, we can select the appropriate $\bD(y)$ to achieve the desired fairness criterion.

      In our example at the beginning of this subsection, there were 2 boundaries of interest: $v_{12} = [1,1]$ and $v_{23} = [1,2]$ (see Figure \ref{fig:N3}).  The boundary which matters depends on the system load, $\rho$, as well as the fairness direction, $\theta$.  For $\theta = [2/3,1/3]$, if $\rho$ is in the lower region, $R_{12}$, then the boundary vector of interest is $v_{12}$, between cones $C_1$ and $C_2$.  If $\rho \in R_{23}$, then the boundary vector of interest is $v_{23}$, between cones $C_2$ and $C_3$.   If $\rho \in R_1$ of $\rho \in R_3$, the fairness direction is infeasible.  These regions can be determined for each subset of service vectors by solving for the set of $\rho$ which satisfy Condition 2 of Corollary \ref{cor:feasible}. As long as we can determine which region $\rho$ resides, the MaxWeight Matrix $\bD$ can be specified without precise knowledge of $\rho$.
\end{enumerate}

See Section \ref{sec:numerical} for a numerical example of how the state evolution depends on $\bD$ and $\rho$ in these two cases.

\subsection{Other Fairness Criteria}\label{ssec:otherfairness}
Thus far we have only considered fairness in terms of the ratios at which workloads will grow.  As we have discussed, there are many different notions of fairness.  In particular, proportional and max-min fairness are two widely used definitions of fairness.  We now discuss how to extend our framework to these notions of fairness.

\textbf{Max-Min Fairness} In Max-Min fairness, the objective is to maximize the minimum utility.  If utility is a decreasing function of workload, this would correspond to minimizing the maximum workload.  Because our system is unstable, the workload in each queue will grow without bound.  However, we know that the time-normalized workload has a finite limit so that:
\eq
\lim_{t\rightarrow \infty} X(t) \approx \eta t
\en
Therefore, to minimize the maximum workload, this would correspond to the objective:
\eq
\min_{\eta} \max_q \eta_q
\en
Hence, we can select $\eta \in \Psi(\rho,\cS)$ which minimizes the maximum index to achieve Max-Min fairness.

\textbf{Proportional Fairness} In the traditional sense of Proportional fairness, the amount of service resources a queue is allocated corresponds to the proportion of anticipated resource consumption required by the queue, $\rho_q$.  Hence, the fairness criterion would require that the average proportion of service rate to queue $q$, $s_q$,  is proportional  to $\rho_q$:
\eq
s_q = \frac{\rho_q}{\sum_k \rho_k}
\en
This implies that, loosely speaking, $X_q(t) = [X_q(s) + \rho_q (t-s) - s_qS(t-s)]^+$, where $S \leq \sum_q \rho_q$ is the time-average aggregate service rate (if $S > \sum_q \rho_q$ then the system is stablizable).  As $t \rightarrow \infty$, this implies that
\eq
\lim_{t\rightarrow \infty} \frac{X_q(t)}{t} = K\rho_q
\en
for some constant $K< 1$.  Hence, we can set $\eta = K\rho$ and achieve proportional fairness provided that $\eta$ is a feasible fairness criterion.

Generally speaking, we can achieve any notion of fairness which can be characterized by a vector $\eta$ which satisfies the constraints in Corollary \ref{cor:feasible}.  Note that unlike much of the conventional work on fairness, our framework presupposes that the queueing system is operating in an unstable regime.

\section{Numerical Results}\label{sec:numerical}
In this section, we present some numerical results to demonstrate the performance of MaxWeight Scheduling in a potential real system.  We examine how the backlogs grow and approach the fairness criterion.  All of our results are asymptotic results with $t\go\infty$.  We can see through some numerical simulations how large $t$ must be in practice to approach our asymptotic results.

To start we look at a system with two ($Q = 2$) queues and two ($N = 2$) service vectors given by the numerical examples in Section \ref{ssec:robust}.
\eq
S_1 = \left(
  \begin{array}{c}
    4 \\
    0 \\
  \end{array}
\right),
S_2 = \left(
  \begin{array}{c}
    3 \\
    1 \\
  \end{array}
\right),
\rho = \left(
  \begin{array}{c}
    4 \\
    1 \\
  \end{array}
\right),
\theta = \left(
  \begin{array}{c}
    2/3 \\
    1/3 \\
  \end{array}
\right)
\en
With only $N = 2$ service vectors, there is only one boundary vector;  $\rho$ is such that the fairness criterion is feasible (as in Corollary \ref{cor:feasible}), so the MaxWeight matrix $\bD = \bD_1 = \left(
  \begin{array}{cc}
    1& 0 \\
    0 &2\\
  \end{array}
\right)$ will achieve the fairness criterion.
Our load vector $\rho = [4,1]^T \not \in \cP$.  In each time slot, the number of jobs which arrive to queue $1$ is uniformly distributed on $[0,8]$; for queue $2$ it is uniformly distributed on $[0,2]$.  In this case $\alpha_1 = 1/3, \alpha_2 = 2/3$ so that $\eta = (\rho-\alpha_1 S_1 - \alpha_2S_2)^+ = [2/3,1/3]^T$.

We consider how the workload vector grows for various initial conditions:
\eq
X(0) = \left(
  \begin{array}{c}
    0 \\
    0 \\
  \end{array}
\right),X(0) = \left(
  \begin{array}{c}
    60 \\
    0 \\
  \end{array}
\right),X(0) = \left(
  \begin{array}{c}
    0 \\
    20 \\
  \end{array}
\right)
\en

\begin{figure}[ht]
\centering
\subfigure[Different initial conditions.]{
\includegraphics[scale=.45]{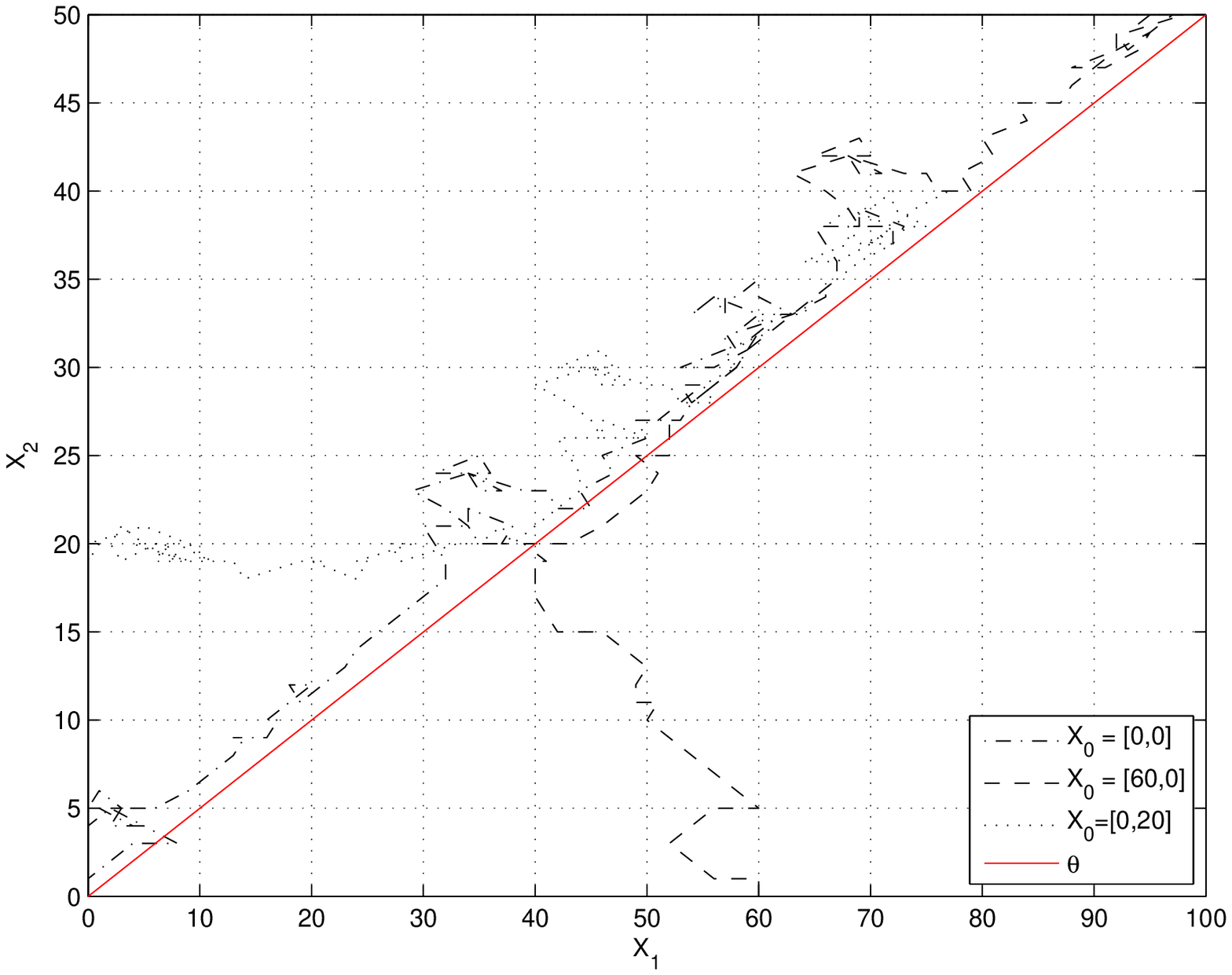}
\label{fig:twoqueues}
}
\subfigure[{Initial condition $X(0) = [0,0]$.}]{
\includegraphics[scale=.45]{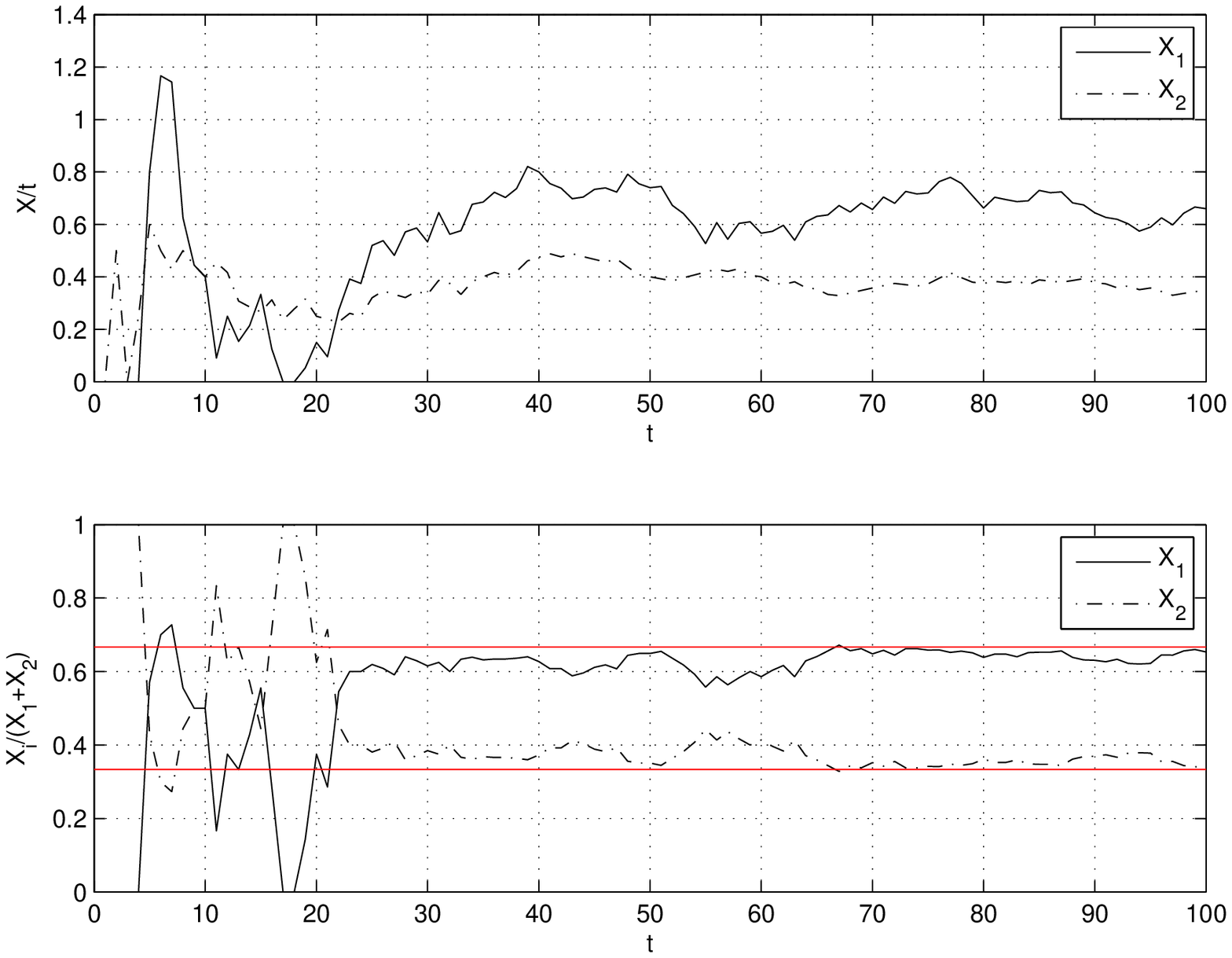}
\label{fig:twoqueues_split}
}
\label{fig:2Q}
\caption{Dyanmics for $N = 2$ service vectors and $Q = 2$ queues.}
\end{figure}

In Figure \ref{fig:twoqueues}, we plot the trajectories of $X(t)$ for the different initial conditions, along with the line $X_1 = 2X_2$.  We can see that all three trajectories converge to the desired fairness vector $\theta_1 = 2/3, \theta_2 = 1/3$.  In Figure \ref{fig:twoqueues_split}, we see the scaled backlogs, $\frac{X_i(t)}{t}$, and the relative backlog, $X_i(t)/\sum_j X_j(t)$, converge starting from initial condition $X(0) = [0,0]$.  Moreover, we see that they quickly achieve the fairness criterion which is shown in red.

\begin{figure}[nh]
\begin{center}
\includegraphics[scale=.7]{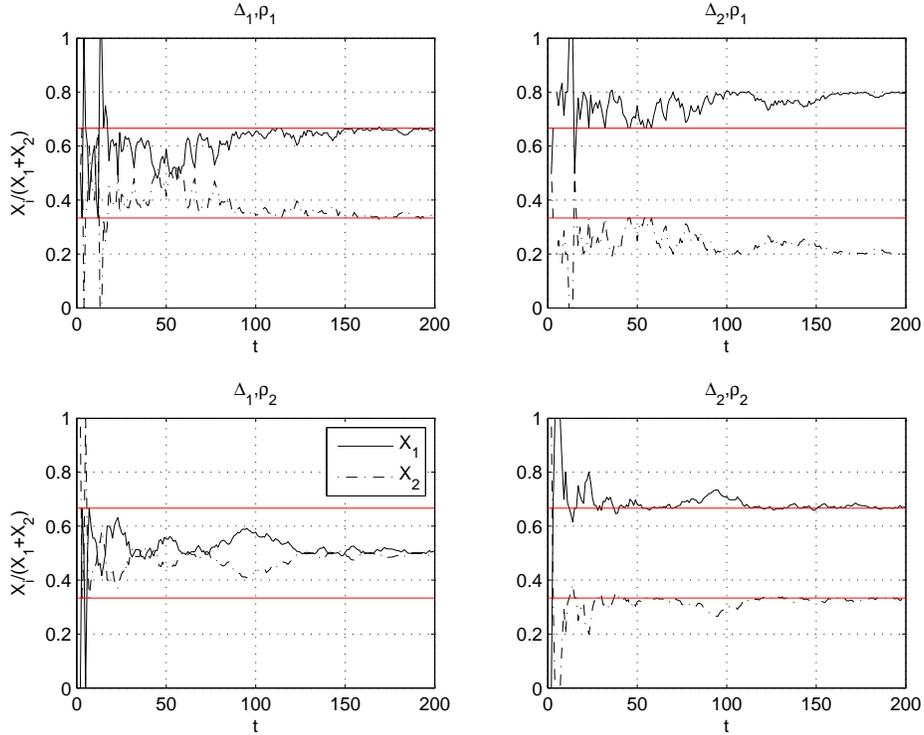}
\caption{Dynamics of 2-queue queueing system with 3 service vectors.} \label{fig:twoqueues3}
\end{center}
\end{figure}

When there are $N = 3$ service vectors, there are $k = 2$ subsets of service vectors and corresponding boundary vectors which satisfy \eqref{eq:relevant}.
\eq
S_1 = \left(
  \begin{array}{c}
    4 \\
    0 \\
  \end{array}
\right),
S_2 = \left(
  \begin{array}{c}
    3 \\
    1 \\
  \end{array}
\right),
S_3 = \left(
  \begin{array}{c}
    0 \\
    4 \\
  \end{array}
\right),
\rho = \left(
  \begin{array}{c}
    4 \\
    1 \\
  \end{array}
\right),
\theta = \left(
  \begin{array}{c}
    2/3 \\
    1/3 \\
  \end{array}
\right)
\en
Now there are two boundary vectors of interest: the one between $C_1$ and $C_2$ as well as the one between $C_2$ and $C_3$.  The boundary which matters depends on $\rho$ as described in Section \ref{ssec:robust}.  From that discussion, we know that   $\bD_1$ and $\bD_2$ are the necessary MaxWeight matrices to move the boundary vectors between $C_1$ and $C_2$ and between $C_2$ and $C_3$, respectively, to fairness direction $\theta$.  We refer to  $\bD(\rho)$ as the necessary MaxWeight Matrix to achieve fairness direction $\theta$, given system load, $\rho$.  Let $\rho = \rho_1 = [4,1]^T$; we use $\bD_1$ as before and we can achieve the desired fairness criterion. On the other hand, if $\rho = \rho_2 = [3,2]^T$, $\bD_1$ will not achieve the desired fairness criterion, but matrix $\bD_2$ will. Hence:
\eq
\bD(\rho_1) = \bD_1 = \left(
  \begin{array}{cc}
    1& 0 \\
    0 &2\\
  \end{array}
\right),
\bD(\rho_2) = \bD_2 = \left(
  \begin{array}{cc}
    1& 0 \\
    0 &4\\
  \end{array}
\right)
\en
 We can see in Figure \ref{fig:twoqueues3} how the asymptotic dynamics of the queues depend on $\rho$ and $\bD$.  The fairness criterion is shown in red.

In the next experiment, we consider a system with $Q=3$ queues and $N=3$ service vectors.  Again there is only one $v \geq 0$ boundary vector to consider.  As long as the fairness direction $\theta$ is feasible there is only one MaxWeight matrix $\bD$ (to a scale factor) to achieve it.  Feasibility depends on $\rho$.
\eq
S_1 = \left(
  \begin{array}{c}
    5 \\
    0 \\
    0 \\
  \end{array}
\right),
S_2 = \left(
  \begin{array}{c}
    0 \\
    5 \\
    0 \\
  \end{array}
\right),
S_3 = \left(
  \begin{array}{c}
    0 \\
    0 \\
    5 \\
  \end{array}
\right)
\en

The desired fairness criterion is
\eq
\theta = \left(
  \begin{array}{c}
    1/2 \\
    1/3 \\
    1/6 \\
  \end{array}
\right)
\en
With some algebra, it is easy to see that as long as the load vector $\rho \not\in \cP$ is such that the fairness direction is feasible, then  the necessary MaxWeight matrix, $\bD$, to achieve it is as follows:
\eq
\bD = \left(
  \begin{array}{ccc}
    2& 0 & 0 \\
    0 &3&0\\
    0 & 0 & 6
  \end{array}
\right)
\en

The system load oscillates between being stable and unstable.  Hence, there are {\em temporary periods of overload}.
\eq
\rho_{\textrm{stable}} = \left(
  \begin{array}{c}
    1 \\
    0 \\
    1 \\
  \end{array}
\right),
\rho_{\textrm{unstable}} = \left(
  \begin{array}{c}
    3 \\
    2 \\
    1 \\
  \end{array}
\right)
\en
The system spends 500 time slots in the stable mode--$\rho = \rho_{\textrm{stable}}$--then switches to spend 500 time slots in the unstable mode--$\rho= \rho_{\textrm{unstable}}$.  Arrivals to queue $q$ in each time slot are uniformly distributed between $[0,2\rho_q]$.  When the system is in the stable mode, MaxWeight Scheduling should stablize the workload.  When it is in the unstable mode, MaxWeight Scheduling should achieve the desired fairness criterion.

\begin{figure}[h]
\begin{center}
\includegraphics[scale=.75]{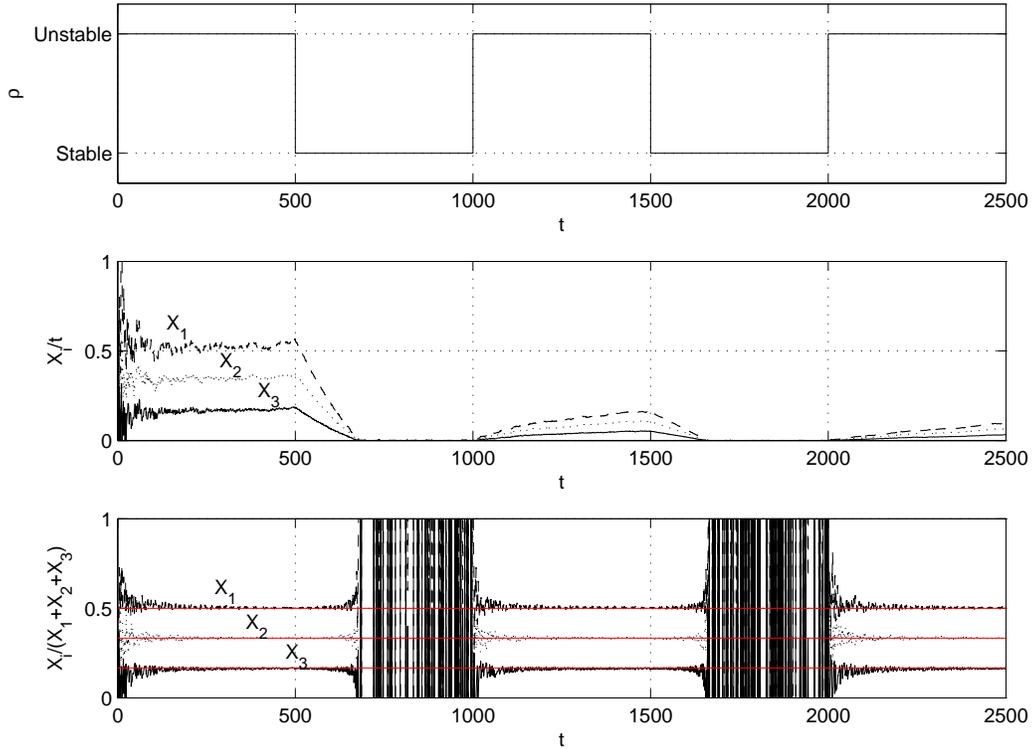}
\caption{Dynamics of 3-queue queueing system oscillating between stable and unstable modes.} \label{fig:multiqueues}
\end{center}
\end{figure}

Figure \ref{fig:multiqueues} plots the scaled workload, $\frac{X(t)}{t}$, under this unstable/stable system. We can see that for the first unstable period ($t\in [0,500]$), the fairness criterion is quickly achieved.  However, on the next unstable period, the scaled backlogs ($X_i(t)/t$) do not appear to stabilize within the $500$ epoch period.  This is because the workload has been stablized during the stable period and we are scaling by the {\em total} time, not just the time starting from when we enter the period of instability.  Hence, it may actually take a very long time before the scaled backlogs converge. One the other hand, the fairness criteria is quickly achieved.  The third subfigure shows how the backlogs  grow relative to each other: $\frac{X_i(t)}{\sum_jX_j(t)}$. We can see in this figure, the fairness criterion is achieved very quickly (plotted in red) during unstable periods.  During stable periods this relative backlog is not very informative since all the backlogs will grow to zero.

While our fairness results for controlling the backlog were asymptotic, these numerical results suggest that the fairness criterion can be achieved very rapidly.

\section{Discussion and Conclusion}\label{sec:conc}
In many real world systems, traffic load is unpredictable and often
bursty in nature.  In any finite window of time, the system may
enter a period of temporary instability where the rate of incoming
jobs is larger than the rate at which jobs can be serviced.  During
such periods of stress, it is natural to want to allocate limited
service resources across various jobs classes in a fair manner.

In this work, we consider how to fairly serve under-provisioned
queues. We focus on MaxWeight Scheduling policies because they are
simple to implement and behave well during stable periods.  In
particular, MaxWeight policies guarantee finite backlogs when the
system load is within the stability region.   We find that, whenever
the system is overloaded, the backlog approaches a straight line as
the time window during which the system is overloaded increases.
This straight line can be characterized as a fixed point, or
equivalently, as the solution to a simple convex program. As such,
it is straightforward to identify this line, as a function of the
system parameters, the load vector $\rho$, and the MaxWeight matrix
$\bD$.

Considering a fair allocation of resources in overloaded systems to
be such that the backlog grows along a certain {\it direction}
$\theta$, we show how to choose the MaxWeight matrix $\bD$ to
guarantee that the backlog will indeed grow along this direction, as
long as this direction is feasible. Additionally, MaxWeight is shown
to asymptotically minimize the backlog, among all other policies
that achieve the same fairness criterion. As it turns out the choice
of $\bD$ is robust with respect to the load vector $\rho$ in the
sense that for every direction $\theta$, there exists a {\it
partition} of the instability region into subsets, such that the
same matrix $\bD$ will work for all load vectors in the same
subset. In particular, it is sufficient to know what subset $\rho$
belongs to, and it is unnecessary to know the exact value of $\rho$.
Via numerical simulations we see that MaxWeight scheduling policy
performs as expected; it achieves the fairness criterion during
periods of overload and it stabilizes the queues otherwise.

Our analysis relies on geometrical arguments that are applied
directly to individual backlog traces. In particular, no
probabilistic assumptions are made. The only necessary assumption is
that the arriving workload has a well defined long run average. This
approach is also helpful in developing intuition with respect to the
system dynamics.

This work can be extended in various directions. As demonstrated in Section \ref{ssec:infeasfeas}, there exist some fairness criteria that are infeasible via MaxWeight Matching.  So, one might consider whether other policies, such as Projective Cone Scheduling from \cite{ross_ton09}, would expand the set of feasible fairness criteria.  In a similar vein, one might
want to obtain other fairness criteria for unstable systems. Third,
it may be possible to extend this work to networks of parallel
queueing systems, by relying on results from \cite{Shah_Wischik}.
Finally, while we have established convergence of the backlog vector
under vary mild traffic conditions, if more restrictive assumptions
are made (such as Markovian queues) one might be able to obtain
results on the rate of convergence as well.

In this work we take a different view to traditional queueing.
Namely, we focus on the \emph{instability} region.  While it is
certainly desirable to operate systems within the stability region,
there are many real world scenarios where this may not be possible.
Input traffic may surge due to unplanned circumstances.  Service
resources may be reduced due to unavoidable accidents or
catastrophes.  During these periods of temporary instability it is
often necessary to allocate limited resources in a fair manner. Once
the system exits the window of stress, it will be stabilizable and
the natural goals of throughput maximization and cost minimization
can be restored.

\bibliographystyle{ormsv080}
\bibliography{fairness}


\appendix

\bProofof {Lemma \ref{lm:str1}}

We write (using similar arguments like in equations of A.22 to A.27 of \cite{ross_ton09}),
\begin{eqnarray}
\<X(t_c)-X(t'_c),\bD\eta\>& =& \<\sum_{t=t'_c}^{t_c-1}A(t),\bD\eta\>
- \<\sum_{t=t'_c}^{t_c-1}D(t),\bD\eta\> \nonumber \\
&= &\<\sum_{t=t'_c}^{t_c-1}A(t),\bD\eta\> - \sum_{t=t'_c}^{t_c-1}\<D(t),\bD\eta\> \nonumber \\
&=& \<\sum_{t=t'_c}^{t_c-1}A(t),\bD\eta\> - \sum_{t=t'_c}^{t_c-1}
\left[ \sum_{q:\eta_q>0} D_q(t)\bD_{qq}\eta_q + \sum_{q:\eta_q=0} D_q(t)\bD_{qq}\eta_q \right] \nonumber \\
&=& \<\sum_{t=t'_c}^{t_c-1}A(t),\bD\eta\> - \sum_{t=t'_c}^{t_c-1}
\left[ \sum_{q:\eta_q>0} S_q(t)\bD_{qq}\eta_q + \sum_{q:\eta_q=0} S_q(t)\ 0 \right] \nonumber \\
&=& \<\sum_{t=t'_c}^{t_c-1}A(t),\bD\eta\> - \sum_{t=t'_c}^{t_c-1}\<S(t),\bD\eta\>\nonumber \\
 &=& \<\sum_{t=t'_c}^{t_c-1}A(t),\bD\eta\> - \max_{S\in\cS}\<S(t),\bD\eta\> (t_c-t'_c),
\end{eqnarray}
To see the above steps, recall the following. First, $X(t)\in\cV(\eta)$ for every $t\in(t'_c,t_c]$ and any $c$, by assumption. Therefore, $X(t)\in\cK(\eta)$, hence, from \eqref{zaza} we get $D_q(t)=S_q(t)$ for $q\in\cQ$ with $\eta_q>0$, for every $t\in(t'_c,t_c]$ and any $c$. Moreover, $X(t)\in\cC(\eta)$ , hence, from \eqref{tata} we get  $\<S(t),\bD\eta\>=\max_{S\in\cS}\<S,\bD\eta\>$, for all $t\in(t'_c,t_c]$ and any $c$.
\eProof

\bProofof { Lemma \ref{lm:str2}}

Note that $t'_n < t_n$ eventually (for any large $n$), expand the terms as follows:
\eq
\frac{\sum_{t=t'_n}^{t_n-1}A(t)}{t_n-t'_n} =
\frac{\sum_{t=0}^{t_n-1}A(t)}{t_n-t'_n} -
\frac{\sum_{t=0}^{t'_n-1}A(t)}{t_n-t'_n} \\
= \frac{\sum_{t=0}^{t_n-1}A(t)}{t_n}\frac{t_n}{t_n-t'_n} -
\frac{\sum_{t=0}^{t'_n-1}A(t)}{t'_n}\frac{t'_n}{t_n-t'_n} ,
\en
and observe that letting $n\go\infty$ we get
\eq
\lim_{n\go\infty}\frac{\sum_{t=t'_n}^{t_n-1}A(t)}{t_n-t'_n} =\rho \ \frac{1}{\chi} - \rho \ (\frac{1}{\chi}-1) = \rho,
\en
since $\lim_{n\go\infty}\frac{\sum_{t=0}^{T}A(t)}{T}=\rho$. This completes the proof of the lemma. \eProof

\bProofof{Lemma \ref{lm:str3}}
We write
\eq
X(t_m) = \sum_{t=0}^{t_m-1}A(t) - \sum_{t=0}^{t_m-1}D(t)
\en
and observe that $D_q(t)=\min\{S_q(t),X_q(t)\}\le S_q(t)$ for every $q\in\cQ$, hence, $-D_q(t) \ge -S_q(t)$. Therefore, since $\bD$ is diagonal (with positive elements), we have $-\<D(t),\bD\eta\> \ge \<S(t),\bD\eta\> \ge -\max_{S\in\cS}\<S,\bD \eta\>$. Projecting on $\bD\eta$ we get
\eq
\<X(t_m),\bD\eta\> \ge
\<\sum_{t=0}^{t_m-1}A(t),\bD\eta\> -\max_{S\in\cS}\<S,\bD \eta\> (t_m)
\en
Dividing by $t_m$ and letting $m\go\infty$, we get
\eq
\<\mu,\bD\eta\> \ge \<\rho,\bD\eta\> - \max_{S\in\cS}\<S,\bD \eta\> >0
\en
This completes the proof of the lemma.
\eProof

\bProofof{Lemma \ref{lm:str4}}
Indeed (recalling that $\bD$ is positive-definite), we have
\begin{eqnarray}
0\le \<\mu-\eta,\bD(\mu-\eta)\> &=&
\<\mu,\bD\mu\> -2\<\mu,\bD\eta\>+\<\eta,\bD\eta\> \nonumber \\
&\le &\<\mu,\bD\mu\> -2\<\eta,\bD\eta\>+\<\eta,\bD\eta\> =
\<\mu,\bD\mu\> - \<\eta,\bD\eta\>,
\end{eqnarray}
so $\<\mu,\bD\mu\> \ge \<\eta,\bD\eta\>$. But since $\<\eta,\bD\eta\>=\limsup{t\go\infty}\<\frac{X(t)}{t},\bD\frac{X(t)}{t}\>$, we must have $\<\mu,\bD\mu\> = \<\eta,\bD\eta\>$, therefore,  $\<\mu-\eta,\bD(\mu-\eta)\>=0$, which implies $\mu=\eta$. This completes the proof of the lemma.
\eProof

\bProofof{Lemma \ref{lm:str5}}
Rewrite the inequality as $- \left[\< \rho, \bD\eta \> - \max_{S\in\cS} \< S, \bD\eta \> \right] \epsilon + \<\eta,\bD\eta\>  \ge (1-\epsilon) \<\eta,\bD\eta\>$, since $1-\epsilon > 0$. This is equivalent (since $\epsilon >0$) to
\eq
\<\eta,\bD\eta\> \ge \<\rho,\bD\eta\> - \max_{S\in\cS} \< S, \bD\eta \>
\en
But this is true by Lemma \ref{lm:str3} applied to the sequence $\{t_c\}$ with $\lim_{c\go\infty}\frac{X(t_c)}{t_c}=\eta$.
This complete the proof of the lemma.
\eProof

\bProofof{Lemma \ref{lm:eps_bounds}}
We have 2 cases to show since, by definition of $\{s_c\}$ we have that $0\le\frac{t_c-s_c}{t_c}\le 1$.  The first case is A) $\epsilon >0$ and the second is B) $\epsilon < 1$.

{\em A)}
We first show that $\epsilon > 0$. We start by showing that there is no increasing unbounded subsequence $\{t_b\}$ of $\{t_c\}$ such that $\lim_{b\go\infty}\frac{t_b-s_b}{t_b} = 0$, where $s_b = \max \{t_a<t_b\}$.  Note that this also implies that $\lim_{b\go\infty}\frac{s_b}{t_b} = 1$. Arguing by contradiction, suppose it exists. Observe that for every $q \in \cQ$ we have
\eq
-\bar{S}_q (t_b - s_b)
\le
X_q(t_b)-X_q(s_b)
\le
\bar{A}_q (t_b -s_b ),
\en
where $\bar{A}_q<\infty$ is the maximum workload that can arrive in queue $q$ in any time slot (see model in \cite{ross_ton09} for assumption of boundedness) and $\bar{S}_q=\max_{S\in\cS}\{S_q\}<\infty$ is the maximum workload that can be removed from queue $q$ in any time slot. Dividing by $t_b$, letting $b\go\infty$, we get
\eq
\lim_{b\go\infty}\frac{X(t_b)-X(s_b)}{t_b}= 0 = \lim_{b\go\infty}\Bigg [\frac{X(t_b)}{t_b} -\frac{X(s_b)}{s_b}\frac{s_b}{t_b}\Bigg ] = \eta - \psi
\en
which implies that $\eta = \psi$ and establishes the desired contradiction.

{\em B)} We still need to show that $\epsilon\neq 1$. Arguing by contradiction, suppose there exists a subsequence $\{t_i\}$ of $\{t_c\}$ (and corresponding subsequence $\{s_i\}$ of $\{s_c\}$) such that $\lim_{i\go\infty}\frac{t_i-s_i}{t_i}=1$, hence, $\lim_{i\go\infty}\frac{s_i}{t_i}=0$. Applying Lemmas \ref{lm:str1} and \ref{lm:str2} with $\{t'_i\}=\{s_i\}$
\eq
\lim_{i\go\infty}\<\frac{X(t_i)-X(s_i)}{t_i-s_i},\bD\eta\> =
\<\rho,\bD\eta\> - \max_{S\in\cS}\<S(t),\bD\eta\>
\en
It follows that
\begin{eqnarray}\label{aqaq1}
\<\eta,\bD\eta\> &=&\lim_{i\go\infty}\<\frac{X(t_i)}{t_i},\bD\eta\>\nonumber \\
&=&
\lim_{i\go\infty}\<\frac{X(t_i)-X(s_i)}{t_i-s_i}\ \frac{t_i-s_i}{t_i} +
\frac{X(s_i)}{s_i}\ \frac{s_i}{t_i},\bD\eta\> \nonumber\\
&=&
\lim_{i\go\infty}\<\frac{X(t_i)-X(s_i)}{t_i-s_i},\bD\eta\>
\ \frac{t_i-s_i}{t_i} +
\< \frac{X(s_i)}{s_i},\bD\eta\> \ \frac{s_i}{t_i} \nonumber \\
&=&
\left[\<\rho,\bD\eta\> - \max_{S\in\cS}\<S(t),\bD\eta\>\right]\cdot 1
+ \<\psi,\bD\eta\>\cdot 0 \nonumber \\
&=&
\<\rho,\bD\eta\> - \max_{S\in\cS}\<S(t),\bD\eta\>
\end{eqnarray}
Now applying Lemma \ref{lm:str3} on the subsequence $\{s_i\}$ with $\lim_{s_i\go\infty}\frac{X(s_i)}{s_i}=\psi$ we get
\eq
\<\psi,\Delta\eta\> \ge
\<\rho,\bD\eta\> - \max_{S\in\cS}\<S,\bD \eta\> = \<\eta,\bD\eta\>,
\en
using \eqref{aqaq1}. Hence, $\<\psi,\Delta\eta\>\ge\<\eta,\bD\eta\>$, which implies $\psi=\eta$ by Lemma \ref{lm:str4}. But this is impossible since by definition of subsequence $\{s_c\}$, $\psi \not = \eta$, which completes the proof of the lemma.
\eProof

\bProofof{Lemma \ref{fixedpoint}}
Consider a subsequence $\{t_n\}$ such that for each $m$:
\eq
\alpha_m = \lim_{n\go\infty} \frac{\sum_{t=0}^{t_n-1} \mathbf{1}_{\{S(t) = S_m\}}}{t_n}
\en
Note that by definition: $\alpha_m \in [0,1]$ and $\sum_m \alpha_m \leq 1$.  Further, because $\rho \not \in \cP$, there exist $q$ and $T < \infty$, such that $X_q(t) > 0$ for all $t > T$; hence, MaxWeight Scheduling policies will never idle for $t >T$ and $\sum_m \alpha_m = 1$.

We have for $q$ such that $\eta_q > 0$:
\begin{eqnarray}
\eta_q &=& \lim_{n\go\infty} \frac{X_q(t_n)}{t_n} \nonumber \\
& = & \lim_{n\go \infty} \frac{\sum_{t=0}^{t_n-1} \big [A_q(t) - D_q(t)\big ]}{t_n} \nonumber \\
& = & \lim_{n\go \infty} \frac{X_q(t_o)+\sum_{t=t_o}^{t_n-1} \big [A_q(t) - \sum_m S_{m,q} \mathbf{1}_{\{S(t) = S_m\}}\big ]}{t_n} \nonumber \\
& = & \rho_q - \sum_m \alpha_m S_{m,q}
\end{eqnarray}
Where  $t_o < \infty$ such that for all $t > t_o$, $X(t) \in \cV(\eta)$.  It's existence is given by
Proposition \ref{prop:uniquelimit}.
For $q$ such that $\eta_q = 0$, we have:
\begin{eqnarray}
0= \eta_q &=& \lim_{n\go\infty} \frac{X_q(t_n)}{t_n} \nonumber \\
& = & \lim_{n\go \infty} \frac{\sum_{t=0}^{t_n-1} \big [A_q(t) - D_q(t)\big ]}{t_n} \nonumber \\
& = & \lim_{n\go \infty} \frac{X_q(t_o) + \sum_{t=t_o}^{t_n-1} \big [A_q(t) - \sum_m \min\{X_q(t),S_{m,q}\} \mathbf{1}_{\{S(t) = S_m\}}\big ]}{t_n}\nonumber \\
& \geq & \lim_{n\go \infty} \frac{X_q(t_o) + \sum_{t=t_o}^{t_n-1} \big [A_q(t) - \sum_m S_{m,q} \mathbf{1}_{\{S(t) = S_m\}}\big ]}{t_n}\nonumber \\
& = & \rho - \sum_m \alpha_m S_{m,q}
\end{eqnarray}
Which means that  $\rho_q -  \sum_m \alpha_m S_{m,q} \leq 0$ and
\eq
0 = \eta_q = \big [ \rho_q - \sum_m \alpha_m S_{m,q}  \big ]^+,
\en
which gives us that $\eta = \big [ \rho - \sum_m \alpha_m S_{m}  \big ]^+ $.

Finally, we have to show that if $\alpha_m > 0$, then $\eta \in C_{S_m}$.  We have seen that  $\alpha_m$ is the proportion of time that service vector $S_m$ is used under MaxWeight Scheduling once $X(t) \in \cV(\eta)$ for all $t > t_o$. By Proposition \ref{prop:uniquelimit},  we know that $t_o$ exists.  By contradiction, suppose that $\eta \not \in C_{S_m}$.  This implies that there exists $m' \not = m$ such that $\<\eta,\bD S_{m'}\> > \<\eta,\bD S_{m}\>$. Since $\alpha_m > 0$, we must use $S_m$ for some $t > t_o$.  This  contradicts the definition of $\cV(\eta)$, which by \eqref{tata} says that MaxWeight Scheduling would use $S_{m'}$ rather than $S_m$ which would imply that $\alpha_m = 0$.  Hence, if $\alpha_m > 0$, $\eta \in C_{S_m}$.
\eProof

\bProofof{Lemma \ref{lemma:etarho}}
This follows from Lemma \ref{fixedpoint}.  Replacing $\eta =  [ \rho - \sum_{m=1}^N \alpha_m S_m  ]^+$ we have
\begin{eqnarray}
\<\eta,\bD\eta\> &=& \< \Big [ \rho - \sum_{m=1}^N \alpha_m S_m \Big ]^+,\bD\eta\> \nonumber \\
 &=& \sum_{q:\eta_q>0} \Big [ \rho - \sum_{m=1}^N \alpha_m S_m \Big ]_q\bD_{qq}\eta_q +  \sum_{q:\eta_q=0}  \Big [ \rho - \sum_{m=1}^N \alpha_m S_m \Big ]^+_q\bD_{qq}\eta_q\nonumber \\
 &=& \sum_{q:\eta_q>0} \Big [ \rho - \sum_{m=1}^N \alpha_m S_m \Big ]_q\bD_{qq}\eta_q +  \sum_{q:\eta_q=0}  \Big [ \rho - \sum_{m=1}^N \alpha_m S_m \Big ]_q\bD_{qq}0\nonumber \\
 &=& \<\rho - \sum_{m=1}^N \alpha_m S_m ,\bD\eta\> \nonumber \\
 & = & \<\rho,\bD\eta\>   - \sum_{m=1}^N\< \alpha_m S_m ,\bD\eta\> \nonumber \\
&=& \<\rho,\bD\eta\> - \max_{S\in\cS}\<S,\bD\eta\>
\end{eqnarray}
where the last equality follows from the fact that $\eta$ is a fixed point.
\eProof

\bProofof{Lemma \ref{lemma:onefixedpt}}
This result follows from the KKT conditions for optimality of the convex program \eqref{eq:cvxprog}.  Our goal is to show that if $\eta = (\rho-\sum_m\alpha_mS_m)^+$ is such that for all $m$ with $\alpha_m > 0$, we have that  $\<\eta,\bD S_m\> \geq \<\eta,\bD S_k\>$, then it is a solution to the convex program \eqref{eq:cvxprog}. The KKT conditions are necessary and sufficient for optimality if the objective function is differentiable and Slater's constraint is satisfied \cite{boyd_book}.  Our objective function is clearly differentiable.

Slater's condition says that there exists $r,\alpha$ such that $r < \sum_m\alpha_m S_m$, $r < \rho$, $\alpha_m > 0$ and $\sum_m \alpha_m = 1$.  It is clear that $r= 0, \alpha_m = 1/M$ satisfies this condition.
Since the KKT conditions are necessary and sufficient in this case and there is only one solution, it will follow that there is exactly one fixed point if all fixed points satisfy the KKT conditions.

To examine the KKT conditions, we first rewrite the convex program in \eqref{eq:cvxprog} as an equivalent convex program:
\begin{eqnarray}
\min_{r,\alpha} && \<\rho-r,\bD (\rho - r)\> \nonumber \\
s.t. && r \leq \sum_m\alpha_m S_m\nonumber \\
&& r \leq \rho \nonumber \\
 && \alpha_m \geq 0, \forall m  \nonumber \\
&& \sum_m \alpha_m = 1
\end{eqnarray}
Let $r^*, \alpha_m^*$ (and correspondingly $\eta^* = \rho - r^* = [\rho - \sum_m \alpha_m^* S_m]^+$) be the solution to the preceding convex program.  The KKT conditions for optimality say that all the constraints must be satisfied and:
\begin{eqnarray} \label{eq:KKT}
\nabla\<\rho-r,\bD (\rho - r)\> + \nabla \lambda'(r - \sum_m\alpha_m S_m) + \nabla \lambda'' (r-\rho) - \nabla \lambda \alpha + \nabla v (\sum_m \alpha_m - 1) & = & 0\nonumber \\
\lambda'_q (r_q-\sum_m\alpha_m(S_m)_q) & = & 0\nonumber \\
\lambda''_q (r_q-\rho_q) & = & 0\nonumber \\
\lambda_m\alpha_m &=& 0 \nonumber \\
\lambda,\lambda',\lambda'',v &\geq& 0
\end{eqnarray}
We look at the first condition:
\begin{eqnarray}\label{eq:KKT_helper}
&\nabla_{r_q}: & 2\bD_{qq}(\rho-r^*)_q = \lambda'_q + \lambda''_q, \forall q \nonumber\\
&\nabla_{\alpha_m}:& \sum_q \lambda'_q (S_m)_q = v - \lambda_m  \nonumber \\
&& \sum_q (2(\rho-r^*)_q\bD_{qq}-\lambda_q'')(S_m)_q = v - \lambda_m\nonumber \\
&& 2\<\rho-r^*,\bD S_m\> = v +\sum_q\lambda''_q- \lambda_m
\end{eqnarray}

Now we show that for any fixed point, there exists $\lambda,\lambda',\lambda'',v$ which satisfy the KKT condition \eqref{eq:KKT}.  To do this, suppose we are given a fixed point  $\eta = \rho - r = (\rho - \sum_m \alpha_m)^+$ with
\eq
\alpha_m > 0 \implies \<\rho - r,\bD S_m\> \geq \<\rho - r,\bD S_k\>
\en
We will construct non-negative Lagrange multiplies to satisfy the KKT conditions.

Consider $q$ such that $\rho_q > r_q$.  In order to satisfy the third constraint in \eqref{eq:KKT},  $\lambda_q'' = 0$.  Now if $\rho_q \leq r_q$, in order to satisfy the first constraint in \eqref{eq:KKT_helper} we have that $\lambda_q' + \lambda_q'' =0$.  To ensure that the Lagrange multipliers are non-negative, we must have that $\lambda_q'' = 0$ for all $q$.  Subsequently:
\eq
\lambda_q' = 2\bD_{qq}(\rho - r)_q, \forall q
\en
Consider $\alpha_m > 0$.  To satisfy the fourth constraint in \eqref{eq:KKT}, $\lambda_m = 0$.  Now to satisfy the second constraint in \eqref{eq:KKT_helper}:
\eq
0 \leq 2\<\rho-r,\bD S_m\> = v, \forall m\textrm{ such that }\alpha_m > 0
\en
which also satisfies the non-negativity of $v$.
Consider $\alpha_m = 0$ and $\alpha_k > 0$.  By the assumption that $\rho-r$ is a fixed point:
\begin{eqnarray}
\<\rho-r,\bD S_m\> & \leq & \<\rho-r,\bD S_k\>\nonumber \\
&\implies& \frac{v - \lambda_m}{2} \leq \frac{v}{2} \nonumber \\
&\implies& \lambda_m \geq 0
\end{eqnarray}
Hence, we can satisfy the KKT conditions in \eqref{eq:KKT} with non-negative $\lambda, \lambda',\lambda'',v$. Since any fixed point satisfies the necessary and sufficient KKT conditions, all fixed points are solutions to the convex program.  There is only one solution and so there is only one fixed point.
\eProof

\bProofof{Lemma \ref{lemma:bplace}}
We show this via construction.  Recall that $\hat{\bD}$ is diagonal: $\<{v},\hat{\bD}S\> = \sum_q {v}_q \hat{\bD}_{qq} S_q$.  For $v_q > 0$
\eq
\bD_{qq} = \frac{\hat{\bD}_{qq}{v}_q}{\eta_q} \geq 0
\en
where the positivity comes from the fact that each element is positive.  If ${v}_q = 0$,
\eq
\bD_{qq} = 1
\en
or some other arbitrary positive value.

Now, for any $i_j$ ($j \in [1,m]$) and $k$ the following holds:
\begin{eqnarray}
\<{v},\hat{\bD} S_{i_j}\> & \geq &\<v,\hat{\bD} S_k\> \nonumber \\
&\Rightarrow &\sum_q \hat{\bD}_{qq} {v}_q(S_{i_j})_q \geq \sum_q \hat{\bD}_{qq} {v}_q(S_{k})_q \nonumber \\
&\Rightarrow &\sum_q v_q  \hat{\bD}_{qq}(S_{i_j})_q -\sum_q \eta_q\bD_{qq}(S_{i_j})_q - \sum_q \eta_q  \bD_{qq} (S_{k})_q \nonumber \\
 &&\geq\sum_q {v}_q \hat{\bD}_{qq}(S_{k})_q - \sum_q\eta_q \bD_{qq}(S_{i_j})_q  -\sum_q\eta_q\bD_{qq}(S_{k})_q \nonumber \\
&\Rightarrow & \sum_q [(S_{i_j})_q (\hat{\bD}_{qq} {v}_q - \bD_{qq}\eta_q) - (S_{k})_q \bD_{qq}\eta_q)] \nonumber \\
&&\geq \sum_q [(S_{k})_q( \hat{\bD}_{qq} {v}_q -\bD_{qq}\eta_q) - (S_{i_j})_q \bD_{qq}\eta_q)] \nonumber \\
&\Rightarrow&- \sum_q (S_{k})_q \bD_{qq}\eta_q \geq -\sum_q (S_{i_j})_q \bD_{qq}\eta_q\nonumber \\
&\Rightarrow& \<\eta,\bD S_{i_j}\>  \geq \<\eta,\bD S_{k}\>
\end{eqnarray}
\eProof

\bProofof{ Proposition \ref{prop:control}}
Assume that Conditions 1 and 2 are satisfied.  We show this implies there exists a $\bD$ such that $\lim_{t\go\infty}\frac{X(t)}{t} = \eta$. We first consider Condition 2. This says that $v\geq0$ is the boundary of cones $C_m$ where $\alpha_m > 0$ defined by $\hat{\bD} = \bI$.  By Lemma \ref{lemma:bplace}, we can construct a $\bD$ such that for all $\alpha_m > 0$:
\eq
\<\eta,\bD S_m\> \geq \<\eta,\bD S_k\>
\en
Now, in conjunction with Condition 1, we have that  $\eta\in\Psi(\rho,\cS)$ is a fixed point. By Theorem \ref{th:limit}, we have that
\eq
 \lim_{t\go\infty}\frac{X(t)}{t} = \eta
\en

Now suppose there exists a $\bD$ such that $\lim_{t\go\infty}\frac{X(t)}{t} = \eta$.  By Theorem \ref{th:limit}, $\eta$ is a (the only) fixed point and, hence, satisfies Condition 1.  Now, we show that Condition 2 is satisfied by constructing the necessary $v\geq 0$.  Similar to the construction of $\bD$ in the proof of Lemma \ref{lemma:bplace} we can determine $v$--the boundary induced when $\hat{\bD} = \bI$.  That is:
\eq
v_q = \bD_{qq}\eta_q
\en
which equals $0$ if and only if $\eta_q = 0$.  This $v_q$ satisfies Condition 2.
\eProof

\bProofof{Theorem \ref{th:perf}}
The proof is via contradiction.  Since both MaxWeight Scheduling and this alternative algorithm achieve the fairness criterion, there exists some $\kappa$ such that $\bar{\eta} = \kappa\eta$. Let's assume that $\eta$ is \emph{not} the smallest scaled workload that achieves the fairness criterion so that  $\kappa < 1$.

Since
\eq
\lim_{t\go\infty}\frac{\bar{X}(t)}{t} = \kappa\eta
\en
for any $\epsilon > 0$, there exists $T$ such that for all $t > T$:
\eq
\Bigg |\kappa\eta_q - \frac{\bar{X}_q(t)}{t}\Bigg | < \epsilon
\en
For any $q$ such that $\eta_q >0$, there exists some $T_q <\infty$ such that for all $t > T_q$, $\frac{\bar{X}_q(t)}{t} > 0$ and so after time $T_q$, $\bar{D}_q =\bar{S}_q(t)$, where $\bar{D}(t)$ and $\bar{S}(t)$ are the departures and service vector used in time slot $t$ as defined  in Section \ref{sec:model}. Given $\epsilon > 0$, there exists $T$ such that:
\eq\label{eq:alphaclose}
0 \leq\frac{\sum_{k=0}^{T-1} \bar{S}_q(k)}{T} - \frac{\sum_{k=0}^{T-1} \bar{D}_q(k)}{T} =  \frac{\sum_{k=0}^{T_q-1} \bar{S}_q(k)}{T} - \frac{\sum_{k=0}^{T_q-1} \bar{D}_q(k)}{T} < \epsilon
\en
Define $\bar{\alpha}_m \geq 0$ be defined as:
\eq \label{eq:alphalim}
\bar{\alpha}_m = \lim_{T\rightarrow \infty}\frac{\sum_{t=0}^{T-1}\bf{1}_{\{\bar{S}(t) = S_m\}}}{T}
\en
Note that $\sum_m \bar{\alpha}_m \leq 1$.  $\bar{\alpha}_m$ is the proportion of time service vector $S_m$ is used up to time $t$. If the limit in \eqref{eq:alphalim} does not exists, we can take a sub-limit which is well-defined since $\bar{\alpha}_m \in [0,1]$.

We have for each $q$ such that $\eta_q > 0$:
\begin{eqnarray}\label{eq:closekappa}
\Bigg | \kappa \eta_q - (\rho - \sum_m \bar{\alpha}_m\bar{S}_m)^+_q \Bigg | & \leq &\Bigg | \kappa \eta_q - \rho_q + (\sum_m \bar{\alpha}_m{S}_m)_q \Bigg |\nonumber \\
& = &
\Bigg | \kappa \eta_q - \frac{\bar{X}_q(T)}{T} + \frac{\bar{X}_q(T)}{T} - \rho_q + \Big(\sum_m \bar{\alpha}_m{S}_m\Big)_q \Bigg |\nonumber \\
&\leq& \Bigg | \kappa \eta_q - \frac{\bar{X}_q(T)}{T}\Bigg | + \Bigg |\frac{\bar{X}_q(T)}{T} -  \rho_q + \Big(\sum_m \bar{\alpha}_m{S}_m\Big)_q \Bigg |\nonumber \\
 &<& \epsilon + \Bigg |\frac{\sum_{k=0}^{T-1} A_q(k)}{T} - \frac{\sum_{k=0}^{T-1} \bar{D}_q(k)}{T} -  \rho_q + \Big (\sum_m \bar{\alpha}_m{S}_m\Big)_q \Bigg | \nonumber \\
 &<& \epsilon +\epsilon + \Bigg |\frac{\sum_{k=0}^{T-1} \bar{D}_q(k)}{T}- \frac{\sum_{k=0}^{T-1} \bar{S}_q(k)}{T}+\frac{\sum_{k=0}^{T-1} \bar{S}_q(k)}{T}-\Big(\sum_m \bar{\alpha}_m{S}_m\Big)_q  \Bigg | \nonumber \\
 &<& 2\epsilon + \Bigg |\frac{\sum_{k=0}^{T-1} \bar{D}_q(k)}{T}- \frac{\sum_{k=0}^{T-1} \bar{S}_q(k)}{T}\Bigg|+\Bigg|\frac{\sum_{k=0}^{T-1} \bar{S}_q(k)}{T}-\Big(\sum_m \bar{\alpha}_m{S}_m\Big)_q  \Bigg | \nonumber \\
 &<& 3\epsilon + \Bigg |\frac{\sum_{k=0}^{T-1} \bar{S}_q(k)}{T}- \Big(\sum_m \bar{\alpha}_m{S}_m\Big)_q  \Bigg | \nonumber \\
 & < & 4\epsilon \implies \kappa\eta_q = \rho_q-\sum_m\bar{\alpha}_mS_m
\end{eqnarray}
The second inequality comes from the fact that $\kappa \eta$ is the limit of $\frac{\bar{X}(t)}{t}$.  The third inequality comes from the definition of $\rho$ as the long-term traffic load. The fourth inequality comes from \eqref{eq:alphalim}.  The last equality comes from \eqref{eq:alphaclose}.

Let $\psi =  (\rho - \sum \bar{\alpha}_m S_m)^+ $:
\begin{eqnarray}
\<\kappa\eta,\bD\kappa \eta\> &=& \<\kappa\eta + \psi - \psi, \bD(\kappa\eta  + \psi - \psi)\>\nonumber \\
& =& \<\kappa\eta  - \psi, \bD(\kappa\eta  - \psi)\> + 2\<\kappa\eta  - \psi,\bD\psi\> + \<\psi,\bD\psi\> \nonumber \\
& \geq & 2\<\kappa\eta  - \psi,\bD\psi\> + \<\psi,\bD\psi\> \nonumber \\
& \geq & 2\<\kappa\eta  - \psi,\bD\psi\> + \<\eta,\bD\eta\> \nonumber \\
& = & 2\sum_{q:\psi >0}(\kappa\eta_q  - \psi_q)\bD_{qq}\psi_q +2\sum_{q:\psi =0}(\kappa\eta_q  - \psi_q)\bD_{qq}\psi_q + \<\eta,\bD\eta\> \nonumber \\
& = & 2\sum_{q:\psi >0}0\bD_{qq}\psi_q +2\sum_{q:\psi =0}(\kappa\eta_q  - \psi_q)\bD_{qq}0 + \<\eta,\bD\eta\> \nonumber \\
& = & \<\eta,\bD \eta\>\nonumber \\
&\implies& \<\kappa\eta,\bD\kappa \eta\> > \<\eta,\bD \eta\>\nonumber \\
&\implies& \kappa\geq 1
\end{eqnarray}
where the second inequality comes from Proposition \ref{prop:etaunique} and the fourth equality comes from\eqref{eq:closekappa}.   We have a contradiction to the assumption that $\kappa <1$
\eProof

\end{document}